\documentclass[12pt]{article}

\usepackage{amssymb,amsmath,bm}
\usepackage{mathrsfs}
\usepackage{constants,color}
\usepackage{enumerate}
\usepackage{appendix}
\usepackage{tikz}
\usepackage{ulem}
\def\qed{\hfill \mbox{\rule{0.5em}{0.5em}}}

\newcommand{\be}{\begin{equation}}
\newcommand{\ee}{\end{equation}}
\newcommand{\bes}{\begin{equation*}}
\newcommand{\ees}{\end{equation*}}
\newcommand{\ba}{\begin{aligned}}
\newcommand{\ea}{\end{aligned}}
\newcommand{\bi}{\begin{itemize}}
\newcommand{\ei}{\end{itemize}}

\newcommand{\NN}{\mathbb{N}}
\newcommand{\EE}{\mathbb{E}}
\newcommand{\PP}{\mathbb{P}}

\newcommand{\AC}{\mathcal{AC}}

\newcommand{\clf}{\mathcal{F}}
\newcommand{\clc}{\mathcal{C}}
\newcommand{\cls}{\mathcal{S}}

\newcommand{\clm}{\mathcal{M}}

\newcommand{\RR}{\mathbb{R}}

\newcommand{\p}{\mathbf{p}}
\newcommand{\subsoln}{\bar{T}}
\newcommand{\norm}{T}
\allowdisplaybreaks

\newcommand{\card}{\textit{card}}

\numberwithin{equation}{section}
\newtheorem{assumption}{Assumption}[section]
\newtheorem{theorem}{Theorem}[section]

\newtheorem{proposition}{Proposition}[section]
\newtheorem{lemma}{Lemma}[section]
\newtheorem{definition}{Definition}[section]

\newtheorem{remark}{Remark}[section]

\newtheorem{condition}{Condition}[section]

\newcommand{\initial}{\mathcal{C}}
\newcommand{\cla}{\mathcal{A}}
\newcommand{\clb}{\mathcal{B}}
\title{Splitting Algorithms for Rare Events of Semimartingale Reflecting Brownian motions}
\author{Kevin Leder, Xin Liu, and Zicheng Wang}
\begin{document}
\maketitle
\begin{abstract}
We study rare event simulations of semimartingale reflecting Brownian motions (SRBMs) in an orthant. The rare event of interest is that a $d$-dimensional positive recurrent SRBM enters the set $B_n = \{z\in\RR^d: \sum_{k=1}^d z_k = n\}$ before reaching a small neighborhood of the origin as $n\to\infty$. We show that under a proper scaling and some regularity conditions, the probability of interest satisfies a large deviation principle. We then construct a subsolution to the variational problem for our rare event, and based on this subsolution construct particle based simulation algorithms to estimate the probability of the rare event. It is shown that the proposed algorithm is stable and theoretically superior to standard Monte Carlo for a broad class of positive recurrent SRBMs. 
\end{abstract}
\tableofcontents

\section{Introduction}
In many scenarios one is interested in the dynamics of an entity that is constrained to a subset of its possible state space by some mechanism. For example, when studying population levels of many interacting species (e.g., chemical compounds or biological populations), a natural constraint is that all population levels are non-negative. Alternatively, one may simply be interested in the dynamics of particles in a container that they cannot escape from. The so-called reflecting diffusions are a group of stochastic processes that have been developed to study the phenomena of constrained stochastic dynamics. Roughly speaking, a reflecting diffusion process behaves as an ordinary diffusion process in the interior of its domain, and it is pushed into the domain whenever the boundary of the domain is hit. 

A substantial amount of research on reflecting diffusions has been carried out going back to the seminal works of Skorokhod \cite{skorohod1961stochastic} and Kingman \cite{kingman1962queues} in one dimension. Watanabe \cite{watanabe1971stochastic} and Stroock \& Varadhan \cite{stroock1971diffusion} proved the existence of reflecting diffusions in higher dimensions when the domain is sufficiently smooth. When dealing with a particle constrained to an orthant, one does not have a smooth domain and an important extension is the work of Harrison and Reiman \cite{harrison1981} that studied semimartingale reflecting Brownian motions (SRBMs) constrained to orthants. The SRBMs in orthants arise as diffusion approximations for open queueing networks in heavy traffic; cf. \cite{harrison1987brownian, harrison1993brownian, williams1996approximation}. Some important results on the existence, uniqueness and positive recurrence of SRBMs can be found in the survey paper by Williams \cite{williams1995semimartingale} (see our Section \ref{sec:PF}). In particular, Hobson \& Rogers \cite{hobson1993} and Williams \cite{williams1985recurrence} established the sufficient and necessary condition for the positive recurrence of SRBMs in quadrants. A fluid-based sufficient condition for the positive recurrence of SRBMs in general dimensions was developed in Dupuis and Williams \cite{dupuis1994}, which was shown to be a necessary condition in three dimension by Bramson, Dai, and Harrison \cite{bramson2010positive}.

Since these papers a vast literature on the topic has emerged, and of particular interest is the body of work concerned with the tail probabilities of stationary distributions of SRBMs (cf. \cite{avram2001,dupuis2002time,dai2011reflecting}). Until now all results on tail probabilities have been for either two dimensional systems or special cases of higher dimensional systems. For example, \cite{liang2013} assumed the property of rotational symmetry for three dimensional SRBMs, while \cite{dai2014multi} studied how to ensure the stationary distribution of an SRBM has a product form, and how this can be used to identify the asymptotics of its tail probabilities. The work \cite{el2012} analyzed the properties of the most likely trajectories that give rise to tail events in SRBMs. Given the challenge of developing analytical results for the tail probabilities of SRBMs, our work focuses on developing numerical algorithms to aid in the estimations of these probabilities.

Broadly speaking, there are two types of algorithms for estimating rare event probabilities -- importance sampling and particle methods. The general idea of importance sampling is to simulate the system of interest under a new probability measure that makes the event of interest less rare and then adjust the estimator by the Radon-Nikodym derivative between the original measure and the sampling measure. In particle methods, one simulates multiple copies of the system of interest and uses some selection method to favor copies of the system where the rare event of interest is more likely. In the current work we will focus on a particularly simple particle method called the ``splitting algorithm". In this method, one establishes landmarks in the state space and each time a particle crosses a landmark it is replaced with a random number of offspring. The creation of landmarks where splitting events occur is of course the challenging aspect of creating efficient and stable splitting algorithms. In \cite{dean2009}, Dean and Dupuis showed how one can use a large deviation principle to ensure that the splitting algorithm is stable and efficiently estimates the rare event of interest. There have been several works (e.g., \cite{dean2009, dean2011,blanchet2011analysis}) that have looked at rare event simulations using particle methods, but to our knowledge, our work is the first to study rare event simulations for SRBMs using any method.

In this work we develop splitting algorithms to aid in the estimations of tail probabilities for SRBMs.  We focus on the particular rare event that a positive recurrent SRBM exits the $L^1$ cube of radius $n$ before returning to the cube of radius $\epsilon$. We show that our splitting algorithm can stably estimate this probability for positive recurrent SRBMs and furthermore provides an improvement over standard Monte Carlo method. Given that it is not in general possible to exactly simulate an SRBM, we use the Euler method to simulate an approximate process. We establish conditions on the discretization size to ensure our approximating process satisfies the same large deviation principle as the original SRBM. As a result our splitting algorithm may be used for SRBMs or the approximating processes with identical asymptotic properties. Given the bias of using an approximating process we also analyze the relative error between the tail probabilities under the SRBM versus the approximating process.

The remainder of this work is organized as follows. In Section \ref{sec:PF}, we define the SRBMs, state the basic assumptions, and define the rare event of interest. In Section \ref{subsec:LDA}, we show that the rare event probability satisfies a large deviation principle. The variational problems for SRBMs are studied in Section \ref{sec:vp}, and especially, in Section \ref{VP_3d}, a subsolution to the variational problem in general dimensions is constructed.  In Section \ref{sec:splitting}, we introduce the discretizations of SRBMs, apply the splitting method to the discretized process, and state our main results on the performance of the splitting method for estimating rare events. Finally, in Section \ref{sec:Num}, we present numerical examples.

Throughout the manuscript we shall use the following notation. For a pair of $d$ dimensional vectors $u$ and $v$, we denote their standard inner product by $\langle u,v\rangle=\sum_{i=1}^du_iv_i$ and the standard Euclidean norm by $|u|=\sqrt{\langle u,u\rangle}$. We write $u\le v$ if $u_i\le v_i$ for $i=1, \ldots, d.$ For a positive definite $d\times d$ matrix $A$ and $d$ dimensional vectors $u$ and $v$, define the associated inner product by $\langle u,v\rangle_{A}=\langle u,Av\rangle$, and the associated norm $|u|_A=\sqrt{\langle u,u\rangle_A}.$ The trace of $A$ is denoted by $tr(A)$. For a Polish space $S$, denote by $\mathcal{C}([0,\infty); S)$ the space of continuous functions from $[0,\infty)$ to $S$, $\mathcal{AC}([0, \infty); S)$ the space of absolutely continuous functions from $[0,\infty)$ to $S$, and $\mathcal{C}^2([0, \infty); S)$ the space of twice differentiable functions from $[0,\infty)$ to $S$. For $f\in \mathcal{AC}([0, \infty); S)$, let $\dot f$ denote its gradient vector, and $\ddot f$ its Hessian matrix. For a Borel set $B\subset \RR^d$, and a constant $c\in \RR$, $cB = \{cz\in \RR^d: z \in B\}.$ For a set $K$ of finite number of elements, denote by $\card(K)$ its cardinality. Throughout the paper we will use the following notation for the asymptotic behavior of positive functions:
\begin{align*}
f(t) \sim g(t) & \quad\hbox{if $f(t)/g(t) \to 1$ as $t \to \infty$}, \\
f(t) = o(g(t)) & \quad\hbox{if $f(t)/g(t) \to 0$ as $t \to \infty$}, \\
%f(t) \gg g(t) & \quad\hbox{if $f(t)/g(t) \to \infty$ as $t \to \infty$} \\
f(t) = O(g(t)) & \quad \hbox{if $f(t) \leq C g(t)$ for all $t$}, \\
f(t) = \Theta(g(t)) & \quad \hbox{if $c g(t) \le |f(t)| \leq C g(t)$ for all $t$},
\end{align*}
where $C$ and $c$ are positive constants.

%Let $(\Omega, \mathcal{F}, \mathbb{P}, \{\mathcal{F}_t\}_{t\geq 0})$ be a filtered probability space satisfying the usual conditions. All the random variables and stochastic processes in this section are assumed to be defined on this space.

\section{Problem formulation}\label{sec:PF}

Semimartingale refecting Brownian motions (SRBMs) in the nonnegative orthant were first introduced by Harrison and Reiman \cite{harrison1981}. Roughly speaking, such a diffusion process behaves like a Brownian motion in the interior of the orthant. Upon reaching the boundary, it is reflected back into the interior in some specified direction, so as to constrain the process in the orthant. The rigorous definition of SRBMs is stated as follows.

\begin{definition}[SRBM] \label{SRBM} An SRBM associated with data $(\theta, \Sigma, R)$ in $\RR^d_+$ is a continuous $\{\clf_t\}$-adapted $d$-dimensional process $Z$ defined on a filtered probability space $(\Omega, \clf, \{\clf_t\}_{t\ge 0}, \PP)$ such that 
\begin{align}\label{eqn:srbm}
Z(t) = z + X(t) + \theta t + RY(t),
\end{align}
where $z\in \RR^d_+$, $X$ is a $d$-dimensional driftless Brownian motion with covariance matrix $\Sigma$ and $X(0)=0$, and $Y$ is a $d$-dimensional stochastic process such that for $i=1,\ldots, d$, 
\bi
\item[\rm (i)] $Y_i(0) = 0;$ 
\item[\rm (ii)] $Y_i$ is continuous and non-decreasing; 
\item[\rm (iii)] $Y_i$ can increase only when $Z_i$ is on $F_i = \{x\in \RR^d_+: x_i =0\}$, i.e., $\int_0^\infty Z_i(t) dY_i(t) = 0. $
\ei
\end{definition}

%Let $R_i, i=1, \ldots, d,$ denote the $i$-th column of $R$. For $x\in \partial \RR^d_+$, the set of directions of constraints is defined as 
%\begin{align*}
%r(x) = \left\{\sum_{i=1}^d c_i R_i: \sum_{i=1}^d c_i =1, c_i \ge 0, \ \mbox{and $c_i >0$ only if $x_i=0$} \right\}.
%\end{align*}

The parameters $\theta, \Sigma$, and $R$ are known as the drift, covariance matrix, and reflection matrix of the SRBM.  We assume that $\Sigma$ is invertible, and for ease of notation we let $D=\Sigma^{-1}$. In \cite{taylor1993}, the authors consider the weak solution of the SRBM equation \eqref{eqn:srbm}, that is the existence of the triplet $(X, Y, Z)$ which satisfies \eqref{eqn:srbm}. It is shown in \cite{taylor1993} that there exists an SRBM $Z$ associated with $(\theta, \Sigma, R)$ if and only if $R$ is completely-$\cls$, and in this case, $Z$ is unique in law and is a Feller continuous strong Markov process. A matrix is said to be completely-$\cls$ if for its every $k\times k$ principle submatrix $G$, there exists a $k$-dimensional vector $v_G$ such that $v_G \ge 0$ and $Gv_G >0.$ 
%Throughout the paper, we will make the following assumption. 
%
%\begin{assumption}\label{existence}
%The matrix $R$ is completely-$\mathcal{S}$ and $\Sigma$ is invertible.
%\end{assumption}

%Sufficient conditions on positive recurrence of SRBMs will be introduced in the next section. 
%From \cite{dupuis1994}, suppose that for $\psi\in\clc([0,\infty);\RR^d_+)$ of the form $\psi(t) = z + \theta t$, $z\in \RR^K_+, t\ge 0,$ all the $\phi$ components of the SP solution for $\psi$ satisfies that $\phi(t) \to 0$ as $t\to\infty.$ Then the SRBM $Z$ is positive recurrent, and has a unique stationary distribution.
An important tool in the study of SRBMs is the Skorokhod problem (SP). Let $\mathcal{C}_0([0,\infty); \RR^d)$ $= \{f\in \mathcal{C}([0,\infty); \RR^d): f(0) \in \RR^d_+\}.$
\begin{definition}[SP]
Let $\psi\in \mathcal{C}_0([0,\infty); \RR^d)$. A pair $(\phi, \eta)\in\mathcal{C}([0,\infty); \RR^d_+)\times \mathcal{C}([0,\infty); \RR^d_+)$ is called a solution of the SP for $\psi$ associated with $R$ if the following hold.
\bi
\item[\rm (i)] For all $t\geq 0, \phi(t)=\psi(t)+R\eta(t)\geq 0$.
\item[\rm (ii)] $\eta$ satisfies the following: (a) $\eta_i(0)=0$, (b) $\eta_i(\cdot)$ is nondecreasing, and (c) $\eta_i(\cdot)$ increases only when $\phi_i(\cdot)=0$, that is, 
$\int_0^\infty \phi_i(t)d\eta_i(t) = 0.$
\ei
\end{definition}
Theorem 2.2 in \cite{williams1995} guarantees that there is a solution to the SP for $\psi$ if and only if $R$ is completely-$\cls$. Define 
$$\Gamma: \mathcal{C}_0([0,\infty); \RR^d) \to \mathcal{C}([0,\infty); \RR^d_+); \psi \mapsto \phi.$$ 
If $\Gamma(\psi)$ is unique for each  $\psi\in \mathcal{C}_0([0,\infty), \RR^d)$, then $\Gamma$ is called the {\em Skorokhod map (SM)} associated with the matrix $R$, and we will call $\Gamma(\psi)$ the {\em $R$-regulation} for $\psi$. %Otherwise, if the solution $(\phi, \eta)$ of $\psi$ is not unique, we let $\Gamma(\psi)$ denote the set of all $\phi$ components of the solutions. 
Throughout this work we will make the assumption below to insure that our SRBM satisfies a large deviation principle.
\begin{assumption}\label{continuity}
The SM $\Gamma$ is well defined on $ \mathcal{C}_0([0,\infty); \RR^d)$, and is Lipschitz continuous with respect to the topology of uniform convergence on compact sets.
\end{assumption}
Note that Assumption \ref{continuity} implies that $R$ is completely-$\mathcal{S}$. The work \cite{dupuis1991lipschitz} develops conditions that guarantee the Lipschitz continuity of $\Gamma$.

We require that the SRBM $Z$ is positive recurrent such that the expected time that $Z$ reaches any closed set with positive Lebesgue measure will be finite.
\begin{definition}[Positive recurrence]
The SRBM $Z$ is positive recurrent if for every closed set $F$ in $\RR^d_+$ with positive Lebesgue measure, $\EE_z[\tau_F] < \infty$, where $\tau_F = \inf\{t\ge 0: Z(t) \in F\}$. 
\end{definition}

Based on Theorem 2.6 of \cite{dupuis1994}, we make the following assumption which provides a sufficient condition for the positive recurrence of our SRBM.
\begin{assumption}\label{stability} 
Let $\psi\in\clc_0([0,\infty);\RR^d)$ be of the form $\psi(t) = z + \theta t$, $t\ge 0$. Then $\lim_{t\to\infty}\Gamma(\psi)(t)=0$.
\end{assumption}

For a positive recurrent SRBM $Z$, we introduce a large deviation scaling with the scaling parameter $n\in\mathbb{N}$ as follows: For $n\in \NN,$ define
\bes
Z_n(t) = \frac{Z(nt)}{n}, \ t\ge 0.
\ees
Let $\epsilon\in(0,1)$. Consider two sets $A_n=\{z\in \RR^d_+: \sum_{k=1}^d z_k  \le \epsilon/n\}$ and $B = \{z\in \RR^d_+: \sum_{k=1}^d z_k \geq 1\}$. We are interested in the probability that the process $Z_n$ first enters $B$ before reaching $A_n$ with a starting point $z_n\notin A_n\cup B$. Letting $\tau_n = \inf\{t\geq 0: Z_n(t)\in A_n\cup B\}$, the probability of interest can be defined as
\begin{align}
\label{eq:Prob_of_Interest}
p_n(z_n) = \PP(Z_n(\tau_n) \in B | Z_n(0) = z_n).
\end{align}
The goal of this work is to build algorithms that can stably estimate the probability $p_n(z_n)$ and provide an improvement over standard Monte Carlo. Given the positive recurrence assumption on the SRBM $Z$ it is intuitively obvious that $p_n$ will go to 0 as $n\to\infty$. Furthermore, the exponential decay rate of the probability will clearly be related to the large deviation properties of the sequence $\{Z_n\}$. One minor point that we will have to deal with is that the sets $\{A_n\}_{n\geq 1}$ are dependent on the scaling factor $n$. Ideally we would like to replace the sequence of sets $\{A_n\}_{n\geq 1}$ with the set $\{0\}$; however for positive recurrent SRBMs it is possible that the expected time to hit the point $0$ is infinite. In the next section, we develop the necessary large deviation theory to characterize (in terms of a variational problem) the exponential decay rate of $p_n(z_n)$.

\section{Large deviation principle for the rare events}\label{subsec:LDA}
For $T\ge 0$, using Schilder's theorem and the contraction principle, $\{Z_{n}(t); t\in[0,T]\}$ with initial value $Z_n(0) = z$ satisfies the large deviation principle (LDP) in $\mathcal{C}([0,T]; \RR^d_+)$ with the good rate function 
\begin{align}\label{rate-finite-time}
	I_{T,z}(\phi) &= \inf_{\psi\in \AC_z([0, T];\RR^d),\phi\in \Gamma_T(\psi)} \frac{1}{2} \int_0^T |\dot \psi(t) - \theta|_{D}^2 dt,
	\end{align}
where %the norm $\|x\| = x'\Sigma^{-1}x$ for $x\in \RR^K$, 
$\Gamma_T(\psi)$ is the $R$-regulation for $\{\psi(t); t\in[0,T]\}$, and 
$$
\AC_z([0, T];\RR^d)=\{\psi\in \AC([0, T]; \RR^d): \psi(0)=z\}.
$$
In particular we have the following result for any $T>0$. Its proof is a simple application of the contraction principle and is thus omitted. 

\begin{lemma}\label{SPLDP} Assume that $z_n \to z \in \RR^d_+$ and let $I_{T,z}$ be as in \eqref{rate-finite-time}. Then
\bi 
\item[\rm (i)] for any closed $F\subset C([0,T]; \RR^d_+),$ 
\bes
\limsup_{n\to\infty} \frac{1}{n} \log \PP_{z_n}\left( Z_{n} \in F \right) \le - \inf_{\phi \in F} I_{T,z}(\phi);\ees
\item[\rm (ii)] for any open $G\subset C([0,T]; \RR^d_+),$
\bes
\liminf_{n\to\infty} \frac{1}{n} \log \PP_{z_n}\left( Z_{n} \in G \right) \ge - \inf_{\phi \in G} I_{T,z}(\phi). \ees
\ei
\end{lemma}

However, the probability $p_n(z_n)$ is time independent and we need to extend the previous result to handle infinite time horizon events. A standard technique is to use the stability properties of the stochastic process to establish an upper bound on the size of the stopping time $\tau_n$. The result is present in the following proposition. 
\begin{proposition}\label{hitting} There exists $c\in (0,\infty)$ such that
\[
\limsup_{n\to\infty} \sup_{z\in \RR^d_+} \frac{1}{n} \log \EE_{z/n}\left(e^{c n \tau_n} \right) < \infty.
\]
where $\EE_{z/n}$ denotes the expectation conditioning on $Z_n(0)=z/n$. 
\end{proposition}

We now consider an infinite horizon variational problem which minimizes over all paths that travel from 0 to a given point $z$:
\be\label{vp}
I(z) = \inf_{T\ge 0} \inf_{\psi\in \AC_0([0,\infty);\RR^d),\Gamma(\psi)(T)= z} \frac{1}{2} \int_0^T | \dot \psi(t) - \theta |_{D}^2 dt,
\ee
where $\Gamma(\psi)(T) =z$ if there is a $\phi\in \Gamma(\psi)$ such that $\phi(T) = z$.
%For $z\in\RR^K$ and $\phi\in C([0, \infty), \RR^K)$, we write $|z| = \sum_{i=1}^K z_i$ and $ |\phi|=\sum_{i=1}^K\phi_i$. 
With Proposition \ref{hitting} the finite horizon LDP can be extended to infinite time horizon events.
\begin{theorem}\label{SPLDP} Assume $z_n = z/n\notin A_n \cup B$. Then
\bes
\lim_{n\to\infty} \frac{1}{n} \log p_n(z_n) = -\inf_{x\in B} I(x).
\ees
\end{theorem}

\section{Variational problems}\label{sec:vp}

This section focuses on the study of the variational problem (VP) defined in \eqref{vp}. A triplet $(\psi, \eta, \phi)$ is called a solution to the VP \eqref{vp} if $(\phi, \eta)$ is a solution of the SP for $\psi$, $\phi(T) = z$ for some $T\ge 0$, and 
\[
\frac{1}{2} \int_0^T | \dot \psi(t) - \theta |_{D}^2 dt = I(z). 
\]
The solution $(\psi, \eta, \phi)$ is also referred to as an {\em optimal triple}, and $\psi$ and $\phi$ are referred to as {\em optimal path} and {\em optimal reflected path}, respectively. 
It is known from \cite{avram2001} that the VP in $\RR_+^2$ admits an explicit analytical solution. We summarize this result in Section \ref{VP_2d}. For the VP in $\RR^3_+$, a special case on rotationally symmetric VP is studied in \cite{liang2013}, and general paths properties are studied in \cite{el2012}. In \cite{farlow2013}, the value of the VP in $\RR^d_+$ is characterized as a solution to a partial differential equation. However, the complete analysis of the VP in three or higher dimensional space is still open. In Section \ref{VP_3d}, we apply the path properties from \cite{el2012} to construct a subsolution for the VP in $\RR^d_+$ when $d\ge 3$.

\subsection{Explicit solution for the VP in $\RR^2_+$}\label{VP_2d}

%From previous section, we know that a subsolution to the VP defined by \ref{vp} is required to construct a stable splitting algorithm. In this section, we will present explicit solution for the VP in 2 dimensional orthant. This problem . 
%\\

%The key step to find the explicit solution for the VP in $\RR^2_+$ is to reduce the search space to the space of piecewise linear functions with at most two segments. Then one can analyze the ``locally" optimal paths with a given structure and obtain the globally optimal path by comparison. In an 2 dimensional orthant, there are only three types of paths which start at the origin and terminate at $v$:
%
%\bi
%\item[(i)] A direct triple to $v$ without boundary reflecting.
%\item[(ii)] A two-segment path from the origin to $v$ through face $F_1$, and each segment is linear.
%\item[(iii)] A two-segment path from the origin to $v$ through face $F_2$, and each segment is linear.
%\ei
The VP in $\RR^2_+$ is analytically solved in \cite{avram2001}, in which the solution is characterized by defining ``cones of boundary influence". %It will turn out that the solution to the VP will only depend on which cone $z$ belongs to.
In this section, we summarize the result in \cite{avram2001} and adopt its notation. Without loss of generality, we let 
\begin{equation}\label{R}
R=
\begin{bmatrix}
    1       & r_2  \\
    r_1   & 1
\end{bmatrix}.
\end{equation}
The parameters $\theta$ and $R$ are required to satisfy the following assumption, which is sufficient and necessary for the positive recurrence of $Z$. 

\begin{assumption}\label{pr_2d}
	The drift $\theta$ and reflection matrix $R$ satify that 
	\begin{equation}\label{stability_2d}
	\theta_1+r_2\theta_2^-<0,\ \ \theta_2+r_1\theta_1^-<0.
	\end{equation}
\end{assumption}

For $i=1,2$, we choose a vector $p^i$ which is orthogonal (under the usual Euclidean inner product)
to the $i$th column of $R$, and $| \Sigma p^i |_{D}=1.$
We define the {\em face} $F_i = \{x\in \RR^2_+: x_i =0\},$ and the associated {\em exit velocity} $a^i=\theta-2(\theta'p^i)\Sigma p^i.$ 
The face $F_i$ is said to be {\em reflective} if the $i$th component of $a^i$ is negative, i.e., $a^i_i<0$. Denote $e^1=[0\;\, 1]^T$, $e^2=[1\;\, 0]^T$. Let $n^i$ be a vector that is normal to $F_i$, pointing to the interior of the state space, and normalized such that $|\Sigma n^i|_{D}=1.$ We introduce the key definition from \cite{avram2001} as follows:

\begin{definition}[Cone of influence]\hfill
\begin{itemize}
\item[\rm (i)] When $F_i$ is not reflective, the cone of influence $C_i$ is defined to be an empty set. In this case, the face $F_i$ has no boundary influence on solutions to the VP for any $z\in \RR^2_+$. 
\item[\rm (ii)]When $F_i$ is reflective, the cone of influence $C_i$ is defined to be the cone generated by $e^i$ and $\tilde{a}^i$, 
where 
\begin{equation*}
\tilde{a}^i=\langle a^i,e^i\rangle_{D} e^i-\langle a^i, \Sigma n^i \rangle_{D} \Sigma n^i.
\end{equation*}
That is
\begin{equation*}
C_i=\{ \alpha e^i+\beta \tilde{a}^i: \alpha\ge 0, \beta\ge 0 \}.
\end{equation*}

\end{itemize}
\end{definition}  

For $z\in \RR^2_+,$ the cone of influence $C_i$ is said to be active if $z\in C_i$. We define three cost functions:
\begin{align*}
I^0(z)&=\langle \tilde{a}^0(z)-\theta,  z\rangle_{D},\\
I^1(z)&=\langle \tilde{a}^1-\theta, z\rangle_{D},\\
I^2(z)&=\langle \tilde{a}^2-\theta, z\rangle_{D},
\end{align*}
where $\tilde{a}^0(z)=z|\theta |_{D}/|z|_{D}$. With Assumption \ref{pr_2d}, the authors in \cite{avram2001} showed that an optimal path can always be found among these three paths. The search for such path depends on the location of $z$ and its active cone of influence. The following is the main theorem from \cite{avram2001}. 

\begin{theorem}[Theorem 3.1 of \cite{avram2001}] Under Assumptions \ref{continuity} and \ref{pr_2d}, the following hold. 
	\bi
	\item[\rm (i)] If $z\notin C_1\cup C_2,$ the optimal value of the VP is given by $I^0(z)$.
	\item[\rm (ii)] If $z\in C_1\backslash C_2,$ the optimal value of the VP is given by $I^1(z)$.
	\item[\rm (iii)] If $z\in C_2\backslash C_1,$ the optimal value of the VP is given by $I^2(z)$.
	\item[\rm (iv)] If $z\in C_1\cap C_2,$ the optimal value of the VP is given by $\min \{ I^1(z), I^2(z) \}$.
	\ei
\end{theorem}

\subsection{Subsolution for VP in $\mathbb{R}^d_+$ when  $d\geq 3$}\label{VP_3d}

\noindent In this section we construct a sub-solution to the VP \eqref{vp} in ${\mathbb{R}}^d_+$, $d\ge 3$. We first define $L = \{1,2,...,d\}$. For $J, K\subset L$, define $R_{JK}$ to be a $\card(J)\times \card(K)$ matrix by deleting rows and columns with indices in $L\backslash J$ and $L\backslash K$, respectively. Let $R_j$ denote the $j$th column of $R$. Finally, define the face associated with $K\subset L$ to be $F_K = \{x\in\RR^d_+: x_i =0, i\in K, x_j >0, j\in L \backslash K\}.$

A $d\times d$ matrix is called a nonsingular $\clm$-matrix if it can be expressed as $\mu \mathbb{I}_d - A$, where $\mathbb{I}_d$ is a $d\times d$ identity matrix, $A$ is a $d\times d$ matrix with nonnegative entries, and $\mu$ is greater than the spectral radius of $A$. We also note that $R^{-1} >0$, which in particular implies that $R^{-1}\theta < 0$ under Assumption \ref{m_matrix}.
In this section, we make the following assumption on $R$
\begin{assumption}\label{m_matrix}
The reflection matrix $R$ is a nonsingular $\mathcal{M}-$matrix, and $\theta < 0.$
\end{assumption}
This implies that the SRBM $Z$ is positive recurrent (see \cite{harrison1981reflected}), and also implies Assumptions \ref{continuity} and \ref{stability} (see \cite{chen1996, harrison1981reflected}).

For any $w,v\in \mathbb{R}^d_+, w\neq v,$ define the globally optimal cost from $w$ to $v$ as follows:
\begin{equation}\label{vp_general}
\mathcal{I}(w,v)=\inf_{\psi\in \AC_w([0,\infty);\RR^K)} \inf_{\phi\in\Gamma(\psi), \tau_\phi(v)<\infty} \frac{1}{2} \int_0^{\tau_\phi(v)} |\dot \psi(t) - \theta|_{D}^2 dt,
\end{equation}
where
\begin{align*}
\tau_\phi(v)=\inf\{ t\ge 0:\phi(t)=v \}.
\end{align*}
A path $\psi^*$ which leads to the globally optimal cost is called a globally optimal path from $w$ to $v$. We next define the locally optimal cost between points $w$ and $v$, where we constrain the path to lie on a face $F_K$. In particular we will consider the local cost between points satisfying the following condition with respect to $K\subset L$.

\begin{condition}\label{wvcondition}
	$w\neq v; w_i=v_i=0, \forall i \in K$, and for each $j\notin K$ $w_j\neq 0$ or $v_j\neq 0$. %{\color{red} Is this okay?}
%	there exist no $K'\subseteq I, $ such that $K\subsetneq K', v\in F_{K'}$, and $w\in F_{K'}$. %at least one of $w$ and $v$ is in $F_K$.
\end{condition}

For the pair $v$ and $w$ satisfying Condition \ref{wvcondition} with respect to $K\subset L$, we define the locally optimal cost from $w$ to $v$ with path constrained in $F_K$ as follows:
%\noindent $\mathcal{I}^*_K(w,v)$ is the locally optimal cost with respect to the constrained variational problem where all paths are constrained in $F_K$.
\begin{equation}\label{LOC1}
\mathcal{I}^*_K(w,v)=\inf_{\psi\in \AC_w([0,\infty);\RR^K)} \inf_{\phi\in\Gamma(\psi), \phi(t)\in F_K, \forall t\in (0,\tau_\phi(v)), \tau_\phi(v)<\infty}\frac{1}{2} \int_0^{\tau_\phi(v)} |\dot \psi(t) - \theta|_{D}^2 dt.
\end{equation}
\\

We now introduce the definition of a sub-solution to the VP \eqref{vp}.
\begin{definition}\label{def:sub-solution}
A function $\subsoln: \mathbb{R}^d_+\rightarrow \mathbb{R}$ is called a sub-solution to the VP \eqref{vp} if
\begin{itemize}
\item[\rm (i)]  $\subsoln$ is continuous;
\item[\rm (ii)] $\subsoln(x)\le 0$ for all $x\in B$;
\item[\rm (iii)] $\subsoln(w)-\subsoln(v)\le \mathcal{I}(w,v)$ for all $w,v\in \mathbb{R}^d_+\setminus (\{0\} \cup B).$
\end{itemize}
\end{definition}

The following lemma presents the main idea of constructing a sub-solution. 

\begin{lemma}\label{ss1}
	Let $\subsoln$ be a function that satisfies conditions (i) and (ii) in Definition \ref{def:sub-solution}. If in addition, for all $K\subseteq L$, $\subsoln(w)-\subsoln(v)\le \mathcal{I}^*_K(w,v)$ for all $w,v$ satisfying Condition \ref{wvcondition} with respect to $K$, 
%\begin{itemize}
%	\item[\rm (i)] $\bar{W}$ is continuous,
%	\item[\rm (ii)] $\bar{W}(x)\le 0$ for all $x\in B$,
%	\item[\rm (iii)] for all $K\subseteq I$, $\bar{W}(w)-\bar{W}(v)\le \mathcal{I}^*_K(w,v)$ for all $w,v$ such that $w_i=v_i=0,$ for all $i \in K$,
%\end{itemize}
then it is a sub-solution to the VP \eqref{vp}.
\end{lemma}

%\noindent ElÂ Kharroubi, Ahmed, Abdelhak Yaacoubi, Abdelghani BenÂ Tahar, and Kawtar Bichard (Queueing Systems 70.4 (2012): 299-337) provided an explicit solution to a problem which is very close to \ref{LOC1}: 
%\\

In \cite{el2012}, the authors have investigated a local variational problem which is very similar to the VP \eqref{LOC1}. Specifically, two types of locally optimal paths were investigated: interior escape paths and single segment boundary escape paths. An interior escape path is the one for which no reflection is applied. Let $K$ and $K'$ be two subsets of $L$ such that $K\subset K'$ and $0<\card(K)<\card(K')\le d$.
For $w\in F_{K'}$, $v\in F_K$, an optimal single segment boundary escape path is a path constrained in $F_K$ which has the cost equal to $\mathcal{I}^*_K$. 
The result in \cite{el2012} can be applied to the VP \eqref{LOC1} directly as the proofs are still valid.

In order to state the results in \cite{el2012} on $\mathcal{I}^*_K$, we need some further notation. Given $K\subseteq L$, and two points $v,w\in \mathbb{R}^d_+$ satisfying Condition \ref{wvcondition} with respect to $K$, for each $J\subseteq K$, define
\begin{align*}
B^J & ={(R_{LJ}^T{\Sigma}^{-1}R_{LJ})}^{-1}R_{LJ}^T{\Sigma}^{-1},  \ \ \mbox{for} \ J\neq \emptyset,\\
A^J & = \begin{cases} \mathbb{I}_{d} -R_{LJ}B^J, & \ \mbox{for} \ J\neq \emptyset, \\
\mathbb{I}_{d}, & \ \mbox{for} \ J = \emptyset, \end{cases} \\
\alpha^J(w,v) & =\frac{| A^J\theta |_{D}}{|A^J(v-w)|_{D}}, \\
{\lambda}^J(w,v) & =\frac{| A^J\theta |_{D}}{|A^J(v-w)|_{D}}B^J(v-w)-B^J\theta, \ \ \mbox{for} \ J\neq \emptyset. 
\end{align*}

\begin{theorem}\label{constrain_vp}
Given $K\subseteq L$, and two points $v,w \in \mathbb{R}^d_+$ which satisfy Condition \ref{wvcondition} with respect to $K$, there exists a unique $J_*\subseteq K$ such that  $\langle A^{J_*}R_j, \alpha^{J_*}(w,v)(v-w)-\theta \rangle_{D} \le 0, \  j\in K\backslash J_*$, and when $J_*\neq \emptyset$, each component of ${\lambda}^{J_*}(w,v)$ is positive. Furthermore,
\begin{align*}
\mathcal{I}^*_K(w,v)&=\langle b_K(w,v)-\theta, v-w \rangle_{D}	\\
&=| A^{J_*}\theta |_{D}  | A^{J_*}(v-w) |_{D}-\langle A^{J_*}\theta , A^{J_*}(v-w)\rangle_{D},
\end{align*}
and an optimal path is
\begin{align*}
x_K(t)=w+tb_K(w,v), \ 0 \le t \le \alpha^{J_*}(w,v), %\frac{|A^{J_*}\theta|_{D}}{|A^{J_*} (v-w)|_{D}},
\end{align*}
where \[b_K(w,v)=\begin{cases} \alpha^{J_*}(w,v)(v-w)-R_{LJ_*}\alpha^{J_*}(w,v)B^{J_*}(v-w)+R_{LJ_*}B^{J_*}\theta, & \ \mbox{when} \ J_*\neq \emptyset, \\ \alpha^{J_*}(w,v) (v-w), & \ \mbox{when} \ J_* = \emptyset. \end{cases} \]
\end{theorem}
%{\color{red} What is the proof of this theorem? Is it a specific result in \cite{el2012}?}

From Theorem \ref{constrain_vp}, we can solve the VP \eqref{LOC1} by identifying the set $J_*$ for the pair of points $v$ and $w$. In particular, when $|K|=0$, the problem \eqref{LOC1} only concerns paths in the interior of $\mathbb{R}^d_+$.

\begin{lemma}\label{Kempty}
Given two points $v,w \in \mathbb{R}^d_+$ such that at least one of them is in the interior of $\mathbb{R}^d_+$, we have
\begin{align*}
	\mathcal{I}^*_{\emptyset}(w,v)=| \theta |_{D} | v-w |_{D}-\langle \theta , v-w\rangle_{D},
\end{align*}
and the optimal path is 
\begin{align*}
	x_{\emptyset}(t)=w+tb_{\emptyset}(w,v), \ \ 0\le t \le \frac{|v-w|_{D}}{|\theta|_{D}},
\end{align*}
where $b_{\emptyset}(w,v)=\frac{|\theta |_{D}}{| v-w|_{D}}(v-w)$.
\end{lemma}

%{\color{red} \textbf{The red section here seems to be unnecessary} When $|K|=1$, without loss of generality, we can assume $K=\{1\}$. The following lemma provides a comprehensive condition to determine whether $J=\{1\}$ or $J=\emptyset$. 

%\begin{lemma}\label{Keq1}
%Given $K=\{1\}$, $J=\{1\}$ if and only if $\gamma \langle R_1,v-w\rangle > \langle R_1,\theta \rangle$, where %$\gamma={\| \theta \|}/{\| v-w \|}$.	
%\end{lemma}
%}

When $|K| \ge 1$, we can find $J_*$ by testing all subsets of $K$ through the conditions from Theorem \ref{constrain_vp}.

In the rest of section, we construct a sub-solution for the VP \eqref{LOC1}, which by Lemma \ref{ss1} is sufficient for establishing a subsolution to \eqref{vp}. We start with $\norm(\cdot)$ being a norm on $\RR^d$. It is required that $\norm(\cdot)$ is monotone on $\mathbb{R}^d_+$, that is if $v, w\in \mathbb{R}^d_+$ and $v \le w$, then $T(v)\le T(w)$. (In Section \ref{subsec:splitting}, we let $T(x) = \sum_{i=1}^d |x_i|$ for $x\in\RR^d$.) We will then rescale this norm in an appropriate fashion by the functions $\mathcal{I}_K^*(\cdot,\cdot), K\subset L,$ and shift it to allow for the largest possible value at the origin.

%and consider the maximal ratio between $\norm(v-w)$ and $\mathcal{I}^*_K(w,v)$ for all proper subset $K\subset I$, {\color{red}and all possible $w,v$ such that condition \ref{wvcondition} is satisfied, $|v-w|=1$, and there exists $i\in I\setminus K$ such that $v_i-w_i\ge 0$.} It will be shown that $\mathcal{I}^*_K(w,v)$ is bounded below by some positive number $\epsilon$. Denote by $\mathbf{r}$ the maximal ratio. Then we scale $\norm(\cdot)$ by $\mathbf{r}$ to get $\hat{\norm}(v)={\norm(v)}/{\mathbf{r}}$. At last, we define that  $\subsoln(v)=-\hat{\norm}(v)+\inf\limits_{w\in B}\hat{\norm}(w)$, and prove that it is a sub-solution. 

Consider a nonempty subset $K\subset L$, and two points $v,w$ satisfying Condition \ref{wvcondition}. We note that, from Theorem \ref{constrain_vp}, if $w',v'\in \mathbb{R}^d_+$ also satisfy Condition \ref{wvcondition}, and $v'-w'=\alpha (v-w)$ for some $\alpha\in (0,\infty)$, then $\mathcal{I}^*_K(w',v')=\alpha \mathcal{I}^*_K(w,v)$. Therefore, it is proper to define $\mathcal{I}^*_K(v-w)=\mathcal{I}^*_K(w,v)$ for any $v,w\in \mathbb{R}^d_+$ that satisfy {Condition \ref{wvcondition}}. One note that $\mathcal{I}^*_K(\cdot)$ is defined on $\RR^d.$
%\noindent \textbf{Remark:} For two pairs $w,v$ and $w',v'$ satisfying Condition \ref{wvcondition}, and $v'-w'=\alpha (v-w)$ for $\alpha\in (0,\infty)$, $\mathcal{I}^*_K(w',v')=\alpha \mathcal{I}^*_K(w,v)$. Therefore, it is proper to define $\mathcal{I}^*_K(u)=\mathcal{I}^*_K(w,v)$ for any $v,w$ that satisfy {\color{red} Condition \ref{wvcondition}}.

\begin{lemma}\label{neg}\hfill
	
Given $K \subseteq L$, and two points $v,w \in \mathbb{R}^d_+$ which satisfy Condition \ref{wvcondition} with respect to $K$, we have
\begin{itemize}
\item[\rm (i)] If  $\mathcal{I}^*_K(w,v)=0$, then $J_*=K$, $v_i-w_i <0, \forall i\in L\backslash K$, and $\norm(v)-\norm(w)\le0$.
\item[\rm (ii)] If $\mathcal{I}^*_K(v)=\mathcal{I}^*_K(w)=0$, then $v=cw$ for some $c > 0$.
\end{itemize}
\end{lemma}

%\begin{lemma}\label{lem:neg}
%Whenever $\mathcal{I}^*_K(w,v)=0$, we have $\norm(v)-\norm(w)\le0$.
%\end{lemma} 
%
%\begin{lemma}\label{lem:linear}
%For some proper subset $K\subset I$, and some $v,w$ such that $v_i=w_i=0, i\in K$, if $\mathcal{I}^*_K(v)=\mathcal{I}^*_K(w)=0$, then $v=cw$ for some constant $c$.
%\end{lemma}

For a nonempty subset $K\subset L$, define a direction set $\mathbb{D}_K\subset \RR^d$ such that a vector $v\in \mathbb{D}_K$ if and only if $v$ satisfies the following conditions:
\begin{itemize}
\item[(1)] $v_i=0, \forall i\in K$,

\item[(2)] $\sum_{i=1}^d |v_i|=1$, 
%{\color{red} Zicheng: Can we replace this with $\sum_{i=1}^d |v_i|=1$?. If yes, please make this change.}

\item[(3)] there exists $j\in L\backslash K$ such that $v_j\ge 0$.
\end{itemize}
Notice that $\mathbb{D}_K$ is a closed set and thus compact due to the continuity of the norm $|\cdot |$. By the definition of $\mathcal{I}^*_K(v)$, it is easy to see that $\mathcal{I}^*_K(v)$ is continuous. 
We now show that $\mathcal{I}^*_K$ is bounded from below on the set $\mathbb{D}_K$.
\begin{proposition}\label{prop:pos}
There exists $\varepsilon > 0$ such that for all $K\subset L$ and $v\in \mathbb{D}_K$,  $\mathcal{I}^*_K(v) > \varepsilon$.
\end{proposition}
The above proposition allows us to consider the maximum of the ratio $\norm$ to $\mathcal{I}^*_K$ over the set $\mathbb{D}_K$. Let 
\begin{equation*}
\mathbf{r}=\max\left[ \frac{\norm(v)}{\mathcal{I}^*_K(v)} : \forall K\subsetneq L, \text{ and }\forall v\in \mathbb{D}_K \right],
\end{equation*}
and define our subsolution to be
\begin{align}\label{def_subsoln}
\subsoln(v)=-\frac{\norm(v)}{\mathbf{r}}+\inf\limits_{x\in B}\frac{\norm(x)}{\mathbf{r}}, \ \quad  v\in \RR^d_+. 
\end{align}

\begin{theorem}\label{th:sub-solution}
$\subsoln$ is a subsolution to the VP \eqref{vp}.
\end{theorem}

\section{Splitting algorithms for SRBMs}\label{sec:splitting}
Particle splitting methods and importance sampling methods are the two most widely used methodologies to obtain numerical estimations of probabilities of rare events. Oftentimes, particle splitting methods are used to simulate a class of rare events called first entrance probabilities. The main idea of particle splitting is to partition the state space of a process into a series of nested subsets. Then the rare event of interest can be considered as a nested sequence of events. When a given subset is entered by a simulated particle for the first time, a number of offspring will be generated at the entry point according to the assigned splitting mechanism. All the offspring will follow the same law of the original process independently. More refined versions of splitting algorithms such as RESTART have been introduced in the last decades (see \cite{dean2011} for details). 

In Section \ref{sec:Euler}, we construct the Euler discretization of the SRBM, and show that the discretized SRBM and the SRBM are exponentially equivalent.  In Sections \ref{subsec:splitting} and \ref{subsec:restart}, we use the solution/sub-solution developed in Section \ref{sec:vp} to develop algorithms to estimate the probability for our rare event of interest. Finally in Section \ref{sec:RareEvent} we analyze the performance of our estimation algorithms.

\subsection{Discrete approximation for SRBMs}\label{sec:Euler}
In this subsection, we construct a discrete approximation to the SRBM $Z$. Recall that the $d$-dimensional SRBM $Z$ in $\RR^d_+$ is given by 

\begin{equation*}
Z(t) = z + X(t) + \theta t + R Y(t), \ t\ge 0,
\end{equation*}
where $X$ is a $d$-dimensional driftless Brownian motion with covariance matrix $\Sigma$, and $Y$ is the corresponding reflecting process. Fix $\Delta >0$. Consider the discretized process $ X^\Delta$ and $\tilde X^\Delta$ defined by $ X^\Delta(t) = X((k-1)\Delta)$, and $ \tilde X^\Delta(t) =  X^\Delta(t)+ \theta (k-1)\Delta$ for $t\in[(k-1)\Delta, k\Delta)$, and $k\in\NN$. Define $(Z^\Delta, Y^\Delta)$ as the solution of the Skorokhod problem for $\tilde X^\Delta$ with the initial $z$. Then we have 
\begin{align}\label{discrete_srbm}
Z^\Delta(t) = z +  \tilde X^{\Delta}(t) + R Y^\Delta(t), \ t\ge 0.
\end{align}
%
%Then
%\begin{equation*}
%Z^n(t)=z + X^{\rho^n}(t) + \theta \rho^n(t) + R Y^n(t) \in S, \ t\ge 0,
%\end{equation*}
%is called the Euler projection approximation to $Z$, where $\rho^n(t) = \max\{k/n: k\in \mathbb{N}\cup\{0\}, k/n \le t\}$, $X^{\rho^n}(t) = X_{k/n}$ for $t\in [k/n, (k+1)/n),$ and $Y^n$ is the reflecting process associated with $X^{\rho^n} + \theta \rho^n. $ 
We solve the Skorokhod problem for $\tilde X^\Delta$ by solving a sequence of Linear Complementarity Problems as outlined in \cite{Williamsnotes}.
%Before we present the method to solve (\ref{discrete_srbm}), we introduce the definition of Linear Complementarity Problem. We adapt the notation from \cite{Williamsnotes}.
%
%\begin{definition}
%	(Linear Complementarity Problem). Let $d$ be a positive integer, and $R$ be a $d\times d$ matrix. Fix $u\in R^d.$ A pair $(v,w)\in R^d \times R^d_+$ is a solution of the Linear Complementarity Problem (LCP) for $u$ associated with $R$ if
%	
%	\begin{itemize}
%		\item  $w=u+Rv$,
%		\item  $w'v=0$,
%	\end{itemize}
%	where $w'$ denotes the transpose of $w$.
%\end{definition}
%
%For simplicity, we introduce the following notation: $w=\text{LCP}(u,R)$. Then $Z^\Delta$ can be computed by the following recurrent formula. \cite{Williamsnotes} For $k\in\NN,$
%
%\begin{equation*}
%Z^\Delta(0)=z, \quad Z^\Delta(k\Delta)=\text{LCP}\left(Z^\Delta((k-1)\Delta)+ X^\Delta(k\Delta) - X^\Delta((k-1)\Delta),R\right),
%\end{equation*}
%and 
%\[
%Z^\Delta(t)=Z^\Delta((k-1)\Delta), \quad  t\in [(k-1)\Delta, k\Delta).
%\]
%

Introducing the large deviation scaling to $Z^\Delta$, we define 
\[
Z^\Delta_n(t) = \frac{Z^\Delta(nt)}{n}, \ t\ge 0. 
\]
Let $z_n = z/n$. The first entrance probability associated with $Z^\Delta$ is defined to be 
\[
p^\Delta_n(z_n) = \PP(Z^\Delta_n(\tau^\Delta_n)\in B| Z^\Delta_n(0) = z_n),
\]
where $\tau^\Delta_n = \inf\{t\ge 0: Z^\Delta_n(t) \in A_n \cup B\}.$
In the following, the first proposition establishes the exponential equivalence between $Z_n$ and $Z^\Delta_n$, and the second one characterize the relative error between $p^\Delta_n(z_n)$ and $p_n(z_n)$.
\begin{proposition}\label{exp_mome}
Under Assumption \ref{m_matrix}, if $\Delta \equiv \Delta(n) \to 0$ as $n\to\infty$, then for any $T\ge 0$ and $\varepsilon>0,$
\[
\limsup_{n\to\infty} \frac{1}{n}\log \PP\left(\sup_{0\le s \le T}|Z_n(s) - Z_n^\Delta(s)| > \varepsilon\right) = -\infty, 
\]
and consequently, for $z_n = {z}/{n}\notin A_n \cup B$, 
\bes
\lim_{n\to\infty} \frac{1}{n} \log p^\Delta_n(z_n) = -\inf_{x\in B} I(x).
\ees
\end{proposition}

\begin{proposition}\label{error}  Let $\eta\in(0, 1/2)$ and $\Delta \equiv \Delta(n)$ satisfy $\lim\limits_{n\rightarrow \infty}{n\Delta}^{\eta}=0$, assume that $z_n=z/n\notin A_n\cup B$, and denote $z_n=(z_{n,1},\ldots, z_{n,d})$. Then under Assumption \ref{m_matrix},
\begin{align}\label{error1}
	\left|\frac{\max_{x\in \initial_n}p^\Delta_n(x)}{\max_{x\in \initial_n}p_n(x)} - 1 \right| = O\left(n\Delta^{\eta}\right), 
	\end{align}
where $\initial_n = \{x \in \RR^d_+: \sum_{k=1}^d x_k = \sum_{k=1}^d z_{n,k}\}.$
\end{proposition}
Ideally Proposition \ref{error} would be a statement of the form
\begin{align}
\label{eq:ideal_error}
\left|\frac{p_n^\Delta(z_n)}{p_n(z_n)}-1\right|=O(n\Delta^\eta),
\end{align}
but unfortunately we were not able to prove such a statement. Since $p_n(x)$ will have the same exponential decay rate for all $x\in\mathcal{C}_n$ we see that considering a maximum over the set $\initial_n$ introduces at most a subexponential effect.

\subsection{Particle splitting methods}\label{subsec:splitting}

Following \cite{dean2009}, we first briefly describe the splitting algorithm to estimate the probabilities $\{p^\Delta_n(z_n)\}_{n\geq 1}$. 
Our state space $\mathbb{R}^d_+$ is partitioned according to a nested collection of sets $B= C_0^n\subset C_1^n\subset \cdots \subset C_K^n \cdots$. These sets are often constructed as the level sets of a particular function $V$, which is called the {\em importance function} in \cite{dean2009}, and will be defined in terms of the solution/sub-solution of the VP \eqref{vp}. 

Given the sequence of nested sets $\{C^n_j\}_{j\geq 0}$ and a splitting number $r$, the most simple splitting algorithm works in the following way. 
Define $l^n(x) = \min\{j\ge 0: x \in C^n_j\}$. The algorithm generates a time inhomogeneous branching process with generations indexed by $i\in \{0, 1, \ldots, l^n(z_n)\}.$   
\begin{itemize}
\item For generation $0$, a single particle starts at $z_n$, and it evolves according to the law of $Z_n^\Delta$ with initial condition $Z_n^\Delta(0) = z_n$.
\item Let $N^n_{i-1}$ denote the number of particles in generation $i-1$, and for  $ j=1, \ldots, N^n_{i-1}$, let $X^n_{i-1,j}$  denote the position of the $j$-th particles in generation $i-1$. The $j$-th particle behaves according to the law of $Z_n^\Delta$ with initial condition $Z_n^\Delta(0) = X^n_{i-1,j}$, and it stops moving when reaching either the set $A_n$ or the next level set $C^n_{l^n(X^n_{i-1,j})-1} = C^n_{l^n(z_n) -i}.$
\item If the $j$-th particle in generation $i-1$ reaches $A_n$ first, do nothing. 
\item If the $j$-th particle in generation $i-1$ reaches $C^n_{l^n(z_n)-i}$ first, denoting by $X^n_{i,j}$ its location of entrance to $C^n_{l^n(z_n)-i}$, generate $r$ number of new particles in generation $i$ with position $X^n_{i,j}$. 
\end{itemize}
Once all the generations have been calculated, define the estimator
\begin{align}\label{PA:est}
s_n =r^{-l^n(z_n)}N^n_{l^n(z_n)}. 
\end{align}

Now consider a sub-solution $\subsoln$. If the problem is two dimensional, we let $\subsoln(v)=-I(v)+\inf_{x\in B}I(x)$, where $I(v)$ is the explicit solution to the VP discussed in Section \ref{VP_2d}. It is easy to check that this is indeed a sub-solution, and $\subsoln(0)=\inf_{x\in B}I(x)$. If the problem is three or higher dimensional, we use the sub-solution constructed in \eqref{def_subsoln} of Section \ref{VP_3d}, where the norm $T$ on $\RR^d$ is defined to be $T(x) = \sum_{i=1}^d |x_i|$. Define the importance function $V(z) = \delta \subsoln(z)/\log(r).$ The level sets of the importance function $V$ are defined to be  
\[
L_x = \{y\in \mathbb{R}^d_+: V(y)\le x\}, \ x\ge 0. 
\] 
Now define the sets $C^n_j$ to be the level sets of $V$. More precisely, for some $\delta >0$, let 
\[C^n_j = L_{(j-1)\delta/n}, j \in \NN.\]
Here $\delta$ is referred to as the level of the splitting algorithm. In the proposition below we assume that $z_n=z/n\notin A_n\cup B$.
\begin{proposition}\label{alg_stable}
The splitting algorithm for $Z^\Delta$ is stable, i.e., the total number of particles ever used grows subexponentially as $n\to\infty$. Letting $\{s_n^{(k)}\}_{k\in\NN}$ be an i.i.d sequence of copies of $s_n$, then for each $N\in\NN,$
\begin{align}\label{split_unbiased}
\EE\left[ \frac{1}{N}\sum_{k=1}^N s_n^{(k)}\right] = p_n^\Delta(z_n).
\end{align}
Furthermore, the variance of $s_n$ can be measured using the following rate of decay of the second moment. 
\begin{align}\label{sp_var}
\lim\limits_{n\rightarrow \infty}-\frac{1}{n}\log \EE_{z_n}[s_n^2]&=\inf_{x\in B} I(x) +\subsoln(0).
\end{align}
\end{proposition}

Proposition \ref{alg_stable} is a combination of Proposition 7 and Theorem 8 from \cite{dean2009}. It should be noted that \cite{dean2009} assumes the rare event of interest is governed by a large deviations principle with a specific form that is different from what we consider in the current paper. However, the proofs of Proposition 7 and Theorem 8 do not depend on that assumption and still apply to the current setting.

\subsection{RESTART splitting}\label{subsec:restart}

The particle splitting algorithm presented in Section \ref{subsec:splitting} is very efficient with regard to the mean and second moment of the number of particles generated if a sub-solution to the associated variational problem can be obtained. However, it can still require a lot of computational effort. The source of inefficiency comes from the fact that most of the particles generated will not end up in the rare event set (in our case, set $B$) \cite{dean2011} . For those particles generated near set $B$ by splitting, most of them will still reach $A_n$ before $B$, and it takes much time to simulate those trajectories. 

Here, we briefly introduce the RESTART splitting algorithm we implemented specifically for our problem. More details about the RESTART algorithm for rare event simulations can be found in \cite{dean2011}. We first define a sequence of importance functions $\{V^n\}$ as follows. 

\begin{equation*}
\tilde{V}^n(y)=\Delta \lfloor \frac{nU(z_n)-nU(y)}{\Delta} \rfloor, \ \mbox{and} \ V^n(y)=0\vee \tilde{V}^n(y),
\end{equation*}
where $\Delta=\log(r)$ with $r$ being a fixed positive integer, $z_n$ is the starting point, and $U(\cdot)$ is a subsolution satisfying the conditions in definition \ref{def:sub-solution}. For each $n$, we define a sequence of nested sets $C^n_0 \supset C^n_1\supset \cdots \supset C^n_{J_n}\supset C^n_{J_n + 1}$, where $C^n_0 = \RR^d_+$, and $C^n_{J_n+1}=\emptyset$, such that $y\in C^n_j\setminus C^n_{j+1}$ if and only if $V^n(y)=j\Delta$. For simplicity, we denote $C^n_j\setminus C^n_{j+1}$ by $D_j$. 

The algorithm starts from an initial particle $y_1$ located at $z_n$. Each particle will be assigned a killing threshold $l$ when it is generated. The killing threshold for the initial particle is set to be $0$. A particle is killed whenever it reaches any $D_j$ such that $j<l$. A killed particle will not be simulated further. During the simulation, whenever a particle moves from $D_j$ to $D_k$, $k>j$, in one step, a total number of $r^{k-j}-1$ (excluding the original particle) offspring will be generated. For each integer $j<\alpha\le k$, $(r-1)r^{\alpha-j-1}$ of the offspring will be assigned the killing threshold $\alpha$. A particle is said to be stopped whenever it reaches set $A_n$ or $B$. Notice that splitting can happen at the moment of a particle being stopped. In such case, all the offspring will also be stopped. A stopped particle will not be simulated further. The algorithm terminates when all the particles are stopped or killed. 

Let $\mathbb{I}$ be the index of those particles which reach set $B$. For a single trial of simulation, $p_n(z_n)$ is then approximated by $r_n=\sum_{i\in I}\exp(V^n(y_i(\tau)))$, where $y_i(\tau)$ is the position of $y_i$ when it is stopped. The following proposition summarizes the main results on the RESTART algorithm for $Z^\Delta_n.$
\begin{proposition}\label{restart_stable}\hfill
\begin{itemize}
\item[\rm (i)]The estimator constructed  from the RESTART algorithm is unbiased. 
\item[\rm (ii)]Suppose that $z_n\rightarrow 0$ with each $z_n\notin (A_n\cup B)$. Then 
	\[
\liminf\limits_{n\rightarrow \infty}-\frac{1}{n}\log \EE_{z_n}[r_n^2]\ge \inf\limits_{y\in B}\{ I(y)+(\subsoln(0)-\subsoln(y))\vee 0  \}
	\]

\item[\rm (iii)] Let $w^n$ denotes the sum of the lifetimes of all the particles simulated in a single trial. Then
\[
\lim\limits_{n\rightarrow \infty}\frac{1}{n}\log \EE_{z_n}[w^n]=0.
\]
\end{itemize}
\end{proposition}
%{\color{red} How does the second moment decay rate compare with regular splitting?}

This proposition is a combination of results from \cite{dean2011}. The RESTART algorithm does not necessarily have a faster decay rate of its second moment than standard splitting. However, due to the lower computational effort per replication it often has superior performance to standard splitting. In particular, there might be a slight increase in variance of the estimator due to the use of RESTART instead of splitting, but this increase is often made up for by not simulating particles that are very unlikely to hit the target set.

\subsection{Estimation of Rare Event Probabilities}\label{sec:RareEvent}
Our goal is to estimate the rare event probability $p_n(z_n)$, and we do this with the estimators $s_n$ or $r_n$. We will focus our discussion on $s_n$, but a similar analysis can be carried out for $r_n$. A standard metric for the quality of the estimator is given by the relative expected mean square error
\begin{align}
\label{eq:MSE}
rMSE(s_n)=\frac{1}{p_n(z_n)^2}\EE\left[\left(p_n(z_n)-s_n\right)^2\right]=\frac{\EE\left[\left(s_n-p_n^\Delta(z_n)\right)^2\right]}{p_n(z_n)^2}+\frac{\left(p_n^\Delta(z_n)-p_n(z_n)\right)^2}{{p_n(z_n)^2}}.
\end{align}
The previous display shows that $rMSE(s_n)$ can be decomposed into a relative variance term and a relative bias term. The relative variance term can be rewritten as
$$
\frac{\EE[s_n^2]}{p_n(z_n)^2}-\left(\frac{p_n^\Delta(z_n)}{p_n(z_n)}\right)^2
$$
From Proposition \ref{error} we know that the second term in the previous display grows at most at subexponential rate. From Theorem \ref{SPLDP} and Proposition \ref{alg_stable} we know that
$$
\frac{\EE[s_n^2]}{p_n(z_n)^2}=\exp\left[n\left(\inf_{x\in B}I(x)-\bar{T}(0)+o(1)\right)\right].
$$
As mentioned earlier for the 2-d problem we have that $\bar{T}(0)=\inf_{x\in B}I(x)$ and thus the relative variance will grow subexponentially in $n$. An unbiased estimator whose relative variance grows subexponentially is said to be efficient, i.e., in two dimensions $s_n$ is an efficient estimator for $p_n^\Delta(z_n)$. For three and higher dimensions we have that 
$$
\bar{T}(0)=\frac{1}{\mathbf{r}}=\min\left[\mathcal{I}_K^*(v): K\subset I, v\in\mathbb{D}_K\right],
$$
and therefore the relative variance term grows with exponential rate
$$
\inf_{x\in B}I(x)-\min\left[\mathcal{I}_K^*(v): K\subset I, v\in\mathbb{D}_K\right],
$$
which is guaranteed to be non-negative due to the subsolution property.

As a point of comparison suppose one were to use standard Monte Carlo to estimate $p_n(z_n)$ and denote the estimator by $mc_n(z_n)$. Note that the bias terms would remain the same but the relative variance term would grow with exponential rate of $\inf_{x\in B}I(x)$. Therefore the estimator $s_n$ has theoretically superior mean square error to $mc_n(z_n)$. Another point to mention, is that analysis of $rMSE(s_n)$ leaves out the computational work for generating a single copy of $s_n$. However, Proposition \ref{alg_stable} guarantees that the expected work per copy of $s_n$ grows subexponentially and thus when normalized for expected work per replicate $s_n$ will still be a superior estimator to $mc_n(z_n)$.

Proposition \ref{error} is unfortunately not strong enough to provide precise control over the relative bias term, but we do know that if it does grow it does so in a subexponential fashion. It is our hope that a stronger statement of the form \eqref{eq:ideal_error} can be achieved.

\section{Numerical Examples}\label{sec:Num}

\subsection{Splitting algorithm}

Before we run the splitting algorithm for a specific example, there are some initial steps needed to be done. We first obtain a sub-solution $\subsoln$ to the associated variation problem. If the problem is two dimensional, we let $T(v)=I(v)$, $\forall v\in \RR_+^2$, where $I(v)$ is the explicit solution to the VP discussed in Section \ref{VP_2d}. Then we let $\subsoln(v)=-T(v)+\inf_{x\in B}T(x)$. If the problem is three or higher dimensional, we can find the scaling factor $\mathbf{r}$ by solving the corresponding optimization problem with the help of any optimization solver, and then construct the sub-solution as defined in (\ref{def_subsoln}) of Section \ref{VP_3d}.
%{\color{red} Zicheng: Please be specific here}. 
Then we choose a positive integer $r>1$ as the number of new particles when a splitting occurs, and choose $\delta\in (0,1]$ as the level size. The importance function is then obtained as $V(z)=\delta\bar{T}(z)/\log(r)$. We then define a collection of level sets $\{ C_j^n:j\in \NN \}$, where $C_0^n=B, C^n_j=L_{(j-1)\delta/n}, j\in \NN$, where $L_a=\{ z:V(z)\le a \}$. Finally, for a point $x$, its level can be calculated through the level function $l_n(x)=\min \{ j\ge 0; x\in C^n_j \}$.

%{\color{red}What norm $T$ do you use? I used euclidean norm for the simulation. Xin suggested that we use 1-norm. I will update all the results here after I finish all the simulation in 1-norm} 
We conduct numerical experiments for one 2-dimensional SRBM, and one 3-dimensional SRBM. The first splitting simulation used 1,000 runs with step size $1/(1000*n)$. The second splitting simulation used 200,000 runs with step size $1/(1000*n)$. Note that our step size does not comply with the requirement in Proposition \ref{error} to achieve the asymptotic error bound. However, for every value of $n$ used in the following simulations, $1/(1000*n)\le 1/(n^{2.5})$, where the latter step size satisfies the requirement in Proposition \ref{error}. In this way, we increase the accuracy of simulation results when $n$ is relatively small, and the polynomial factor is not too large as $n$ increases. Simulations are run in Matlab with computer having Intel i7-7500U CPU @ 2.70GHz, and RAM 8.00 GB(7.88 GB useable). 

For a $2$-dimensional SRBM, we consider the initial position $(0.1,0.1)$, $R=[1,0;-1,1]$, $\theta=[-2;1]$, $\Sigma=[1,0;0,1]$, $A_n=\{z\in \RR^2_+: \sum_{k=1}^2 z_k  \le 0.15/n\}$, $B=\{z\in \RR^2_+: \sum_{k=1}^2 z_k  \ge 1\}$.

\begin{center}
	\begin{tabular}{| l | l | l | l |}
		\hline
		n  & 5 & 10 & 15\\ \hline
		Estimate Value & $6.3000\times10^{-4}$ & $5.5500\times10^{-6}$ &  $3.8000\times10^{-8}$ \\ \hline
		Time Taken(second)  & $6.4$ & $386.8$ & $871.3$\\ \hline
		Standard Error  & $1.59\times10^{-4}$ & $6.63\times10^{-7}$ & $7.06\times10^{-9}$\\ \hline
		95\% C.I.  & $[3.18, 9.42]\times10^{-4}$ & $[4.25, 6.85]\times10^{-6}$ & $[2.42, 5.18]\times10^{-8}$\\ \hline
		Average no. particles  & $1.96$ & $27.99$ & $37.94$\\ \hline
		Std no. particles  & $5.12$ & $72.32$ & $132.98$\\ \hline
		Max no. particles  & $55$ & $496$ & $1288$\\ \hline
	\end{tabular}
\end{center}

For a $3$-dimensional SRBM, consider the initial position $(0.1,0.1,0.1)$, $R=[3,-1,-1;$ $-1,2,-1;-1,-1,2]$, $\theta=[-2;-1;-1]$, $\Sigma=[2,1,1;1,2,1;1,1,3]$, $A_n=\{z\in \RR^3_+: \sum_{k=1}^3 z_k  \le 0.15/n\}$, $B=\{z\in \RR^3_+: \sum_{k=1}^3 z_k  \ge 1\}$.

\begin{center}
	\begin{tabular}{| l | l | l | l |}
		\hline
		n  & 5 & 10 & 15\\ \hline
		Estimate Value & $8.5600\times10^{-3}$ & $4.1250\times10^{-4}$ &  $2.1500\times10^{-5}$ \\ \hline
		Time Taken(second)  & $66.18$ & $405.9$ & $117.9$\\ \hline
		Standard Error  & $2.06\times10^{-4}$ & $1.49\times10^{-5}$ & $3.28\times10^{-6}$\\ \hline
		95\% C.I.  & $[8.15, 8.96]\times10^{-3}$ & $[3.83, 4.42]\times10^{-4}$ & $[1.51, 2.79]\times10^{-5}$\\ \hline
		Average no. particles  & $1.00$ & $1.43$ & $1.02$\\ \hline
		Std no. particles  & $0.00$ & $1.92$ & $0.392$\\ \hline
		Max no. particles  & $1$ & $10$ & $10$\\ \hline
	\end{tabular}
\end{center}

\begin{remark}
Splitting algorithm does not always work well for 3-dimensional SRBM due to the fact that the subsolution might not be good enough in some cases. 
\end{remark}

We also ran the standard Monte Carlo simulation for the 3-dimensional example with $n=5$ to compare with the results generated by splitting algorithm. We ran splitting algorithm $100,000$ times so that $100,000*r^{l(x_0)-1}=1,000,000$. We also ran standard Monte Carlo $1,000,000$ times.

\begin{center}
	\begin{tabular}{| l | l | l |}
		\hline
		n=5 & Splitting & Std MC \\ \hline
		Estimate Value & $8.430\times10^{-3}$ & $8.439\times10^{-3}$ \\ \hline
		Time Taken(second)  & $32.52$ & $288.86$ \\ \hline
	\end{tabular}
\end{center}

\subsection{RESTART algorithm}

We applied the RESTART algorithm to simulate the same $2$-dimensional example as in the previous section. As expected, the RESTART algorithm required much less computational effort. Hence, we could run more trials, and simulate the cases for high value of $n$ which would be very time consuming when splitting algorithm is applied.

Here we give numerical results for a $2$-dimensional SRBM. The simulation used 10,000 runs with step size $1/(1000*n)$. Simulations are run in Matlab with computer having Intel i7-7500U CPU @ 2.70GHz, and RAM 8.00 GB(7.88 GB useable).  

For $2$-dimensional SRBM, consider the initial position $(0.1,0.1)$, $R=[1,0;-1,1]$, $\theta=[-2;1]$, $\Sigma=[1,0;0,1]$, $A_n=\{z\in \RR^2_+: \sum_{k=1}^2 z_k  \le 0.15/n\}$, $B=\{z\in \RR^2_+: \sum_{k=1}^2 z_k  \ge 1\}$.

\begin{center}
	\begin{tabular}{| l | l | l | l |}
		\hline
		n  & 5 & 10 & 15\\ \hline
		Estimate Value & $8.853\times10^{-4}$ & $5.133\times10^{-6}$ &  $3.796\times10^{-8}$ \\ \hline
		Time Taken(second)  & $435.8$ & $492.4$ & $660.7$\\ \hline
		Standard Error  & $9.57\times10^{-5}$ & $9.24\times10^{-7}$ & $7.03\times10^{-9}$\\ \hline
		95\% C.I.  & $[6.977, 10.73]\times10^{-4}$ & $[3.322, 6.944]\times10^{-6}$ & $[2.419, 5.174]\times10^{-8}$\\ \hline
		Average no. particles  & $36.55$ & $114.61$ & $256.79$\\ \hline
		Std no. particles  & $224.37$ & $902.52$ & $2.26\times10^{3}$\\ \hline
		Max no. particles  & $4.671\times10^{3}$ & $2.675\times10^{4}$ & $5.949\times10^{4}$\\ \hline
	\end{tabular}
\end{center}

\begin{center}
	\begin{tabular}{| l | l | l | l |}
		\hline
		n  & 20 & 30 & 40\\ \hline
		Estimate Value & $2.611\times10^{-10}$ & $1.132\times10^{-14}$ &  $3.146\times10^{-19}$ \\ \hline
		Time Taken(second)  & $852.4$ & $965.0$ & $1107.0$\\ \hline
		Standard Error  & $5.30\times10^{-11}$ & $3.19\times10^{-15}$ & $1.10\times10^{-19}$\\ \hline
		95\% C.I.  & $[1.573, 3.650]\times10^{-10}$ & $[0.507, 1.757]\times10^{-14}$ & $[0.992, 5.299]\times10^{-19}$\\ \hline
		Average no. particles  & $415.3$ & $669.8$ & $845.31$\\ \hline
		Std no. particles  & $4.06\times10^{3}$ & $8.58\times10^{3}$ & $1.39\times10^{4}$\\ \hline
		Max no. particles  & $1.114\times10^{5}$ & $3.267\times10^{5}$ & $7.602\times10^{5}$\\ \hline
	\end{tabular}
\end{center}

\section{Discussion}
%{\color{red}Discuss what we have learned in this project. For example, RESTART is better than regular splitting. 3D problems are hard etc...
%Possible extensions.}

In this project, we developed particle based simulation algorithms for estimating rare event probability related to SRBMs in a non-negative orthant. We found that splitting-type algorithms provide much better performance than standard Monte Carlo method. Furthermore, we found that the RESTART splitting algorithm is superior to regular splitting algorithm in terms of operating time. Following \cite{dean2009, dean2011} our algorithms are based on subsolutions to the variational problems associated with our rare events. In particular, the work \cite{avram2001} developed solutions to these variational problems in two dimensions. For three and higher dimensions we construct new subsolutions to the variational problem. Due to the results of \cite{dean2009, dean2011} we are able to show that our splitting algorithm has superior performance to standard Monte Carlo in two and higher dimensions.

Future directions for this work include developing subsolutions in three and higher dimension that are closer to the solution, and thus enable simulation algorithms with better performance. Another possible direction is to consider more general reflecting diffusions, i.e., reflecting diffusions with state dependent drift and variance.

\section{Proofs} 
We provide all the proofs in this section. 

\subsection{Proofs for Section \ref{subsec:LDA}}
To show Proposition \ref{hitting}, we introduce the following Lyapunov function introduced in \cite{dupuis1994}. Recall that $R_i$ is the $i$th column of the reflection matrix $R$, and for a $d\times d$ matrix $M$, its trace is denoted as $tr(M)$. 

\begin{theorem}[\cite{dupuis1994}]\label{lyp}
	Under Assumption \ref{stability}, and the assumption that $R$ is completely-$\mathcal{S}$, there exists a Lyapunov function $L: \RR^K_+ \to \RR_+$ such that 
	\begin{itemize}
		\item[\rm (i)] $L \in \clc^2(\RR^d_+\backslash \{0\})$;
		\item[\rm (ii)] Given $M_0>0$, there exists an $M_1>0$ such that when $|z|\ge M_1$, we have $L(z) \ge M_0$;
		\item[\rm (iii)] Given $\varepsilon >0,$ there exists a $M_2 >0$ such that when $|z|\ge M_2$, we have $|\ddot L(z)|\le \varepsilon$;
		\item[\rm (iv)] There exists $c_0 >0$ such that 
		\begin{align*}
		&\langle \dot L(z), \theta \rangle \le - c_0, \ \ \mbox{for all $z\in \RR^d_+$.}\\
		&\langle \dot L(z), R_i \rangle \le - c_0, \ \ \mbox{for all $z\in \partial \RR^d_+\backslash \{0\}$.}
		\end{align*}
		\item[\rm (v)] $L(\alpha z) = \alpha L(z)$ for all $\alpha \ge 0$ and $z\in \RR^d_+.$
	\end{itemize}
\end{theorem}
Some consquences of the above properties are obtained in \cite{budhiraja2007}, which are listed as follows. 
\begin{itemize}
	\item[\rm (vi)] For every $M_3>0$, there exists a $\gamma(M_3)>0$ such that $\sup_{|z|\le M_3} L(z) \le \gamma(M_3).$
	\item[\rm (vii)] There exists an $M_4 >0$ such that $\sup_{z\in\RR^d_+ \backslash \{0\}}| \dot L(z)|\le M_4$.
	\item[\rm (viii)] There exists $0<c_1 \le c_2< \infty$ such that $c_1 |z| \le L(z) \le c_2 |z|.$
\end{itemize}

\begin{lemma}[\cite{budhiraja2007}]\label{martingale-lemma}
	Fix $x\in \RR^d_+$ and $\Delta >0$. For $m\in\NN$, define 
	\begin{align*}
	\nu_m = \sup_{(m-1)\Delta \le t\le m\Delta} \left|\int_{(m-1)\Delta}^t \langle \dot L(Z(s)), d X(s) \rangle \right|.
	\end{align*}
	Then for any $\kappa>0$ and $m, n\in \NN, m\le n$, we have 
	\begin{align*}
	\EE_x \left(\exp\left\{\kappa \sum_{i=m}^n \nu_i\right\} \right) \le \left( 2\sqrt{2} \exp\left\{\kappa^2 M_4^2 \gamma_0 \Delta\right\}\right)^{n-m+1},
	\end{align*}
	where $M_4$ is as in Theorem \ref{lyp} (vii) and $\gamma_0$ is a positive constant only depending on the norm of $\Sigma^{-1}.$
\end{lemma}

{\bf Proof of Proposition \ref{hitting}.}
Define $\tau_\epsilon = \inf\{t\ge 0: \sum_{k=1}^d Z_k(t) \le \epsilon\}$ and then $n\tau_n \le \tau_\epsilon$. Thus it suffices to show that there exists $c>0$ such that
\[
\limsup_{n\to\infty} \sup_{z\in \RR^d_+} \frac{1}{n} \log \EE_{z}\left(e^{c \tau_\epsilon} \right) < \infty.
\]
The proof idea is adapted from that of Theorem 4.3 of \cite{budhiraja2007}. From Theorem \ref{lyp} (iii), choose $M_2>0$ such that when $\sum_{k=1}^d z_k \ge M_2$, $tr(\ddot L(z)\Sigma) < c_0$, where $c_0$ is as in Theorem \ref{lyp} (iv). Let $M_0 > M_2$, and define $\kappa_\epsilon = M_0/\epsilon.$ For $n\in\NN$, consider the set $A_n = \{\omega \in \Omega: \inf_{0\le s \le n\Delta} \sum_{k=1}^d Z_k(s) > \epsilon\}.$ For $\omega \in A_n$ and $m\le n$, from Ito's formula, we have 
\begin{align*}
L(Z(m\Delta)) & = L(Z((m-1)\Delta)) + \int_{(m-1)\Delta}^{m\Delta} \frac{1}{2}tr(\ddot L(Z(s))\Sigma) ds \\
& \quad + \int_{(m-1)\Delta}^{m\Delta} \langle \dot L(Z(s)), \theta \rangle ds + \int_{(m-1)\Delta}^{m\Delta} \langle \dot L(Z(s)), dX(s) \rangle \\
& \quad + \sum_{i=1}^d \int_{(m-1)\Delta}^{m\Delta} \langle \dot L(Z(s)), R_i \rangle dY_i(s). 
\end{align*}
Multiplying both sides of the above equation by $\kappa_\epsilon$, from Theorem \ref{lyp} (iv), (v) and (viii), we have for $m\le n$, on $A_n$, 
\begin{align*}
c_1M_0 <\kappa_\epsilon L(Z(m\Delta)) \le \kappa_\epsilon L(Z(m-1)\Delta) - \frac{c_0\Delta}{2} + \kappa_\epsilon \nu_m,
\end{align*}
which implies that
\begin{align}\label{An}
c_1M_0 <\kappa_\epsilon L(Z(n\Delta)) \le \kappa_\epsilon L(z) - \frac{c_0n\Delta}{2} + \kappa_\epsilon\sum_{m=1}^n\nu_m.
\end{align}
From \eqref{An} and using Lemma \ref{martingale-lemma}, we have for any $\kappa > 0,$
\begin{align*}
\PP_z(A_n) & \le \PP_z\left(\kappa_\epsilon\sum_{m=1}^n\nu_m\ge c_1 M_0 + \frac{c_0n\Delta}{2}  - \kappa_\epsilon L(z) \right) \\
&  = \PP_z\left(\kappa \kappa_\epsilon\sum_{m=1}^n\nu_m\ge \kappa \left(c_1 M_0 + \frac{c_0n\Delta}{2}  - \kappa_\epsilon L(z)\right) \right)  \\
& \le \exp\left\{-\kappa \left(c_1 M_0 + \frac{c_0n\Delta}{2}  - \kappa_\epsilon L(z)\right) \right\} \EE_x\left[\exp\left\{\kappa \kappa_\epsilon \sum_{m=1}^n\nu_m\right\} \right]\\
& \le \exp\{\kappa\kappa_\epsilon L(z)\}\exp\left\{n\Delta[\log(2\sqrt{2})/\Delta + \kappa^2 \kappa_\epsilon^2M_4^2 \gamma_0 - \kappa c_0/2] - \kappa c_1 M_0 \right\}. 
\end{align*}
We can choose $\Delta$ large enough and $\kappa$ small enough such that $-\eta \equiv \log(2\sqrt{2})/\Delta + \kappa^2 M_4^2 \gamma_0 - \kappa c_0/2 < 0.$  Therefore, we have 
\begin{align*}
\PP_z(A_n) & \le  \exp\{\kappa\kappa_\epsilon L(z)\} \exp\{-\eta n\Delta\}.
\end{align*}
For $t>0$, there exists $n_0\in\NN$ such that $t\in [n_0\Delta, (n_0+1)\Delta]$, and we have 
\begin{align*}
\PP_z(\tau_\epsilon>t) \le \PP_z(A_{n_0}) & \le  \exp\{\kappa\kappa_\epsilon L(z)\} \exp\{-\eta n_0\Delta\} \le \exp\{\kappa\kappa_\epsilon L(z)+\eta \Delta\} \exp\{-\eta t\}.
\end{align*}
When $z\in \RR^d_+$ satisfies that $\sum_{k=1}^d z_k\le \epsilon$ or $\sum_{k=1}^d z_k \ge n$, we have $\tau_n =0$. When $\epsilon < \sum_{k=1}^d z_k < n$, we have for $c \in (0, \eta)$, 
\begin{align*}
\EE_{z/n}(e^{cn\tau_n}) & \le \EE_z(e^{c\tau_\epsilon})  = 1 + \int_0^\infty c e^{ct}\PP_z(\tau_\epsilon > t) dt \\
& \le 1 +  \frac{c}{\eta - c} \exp\{\kappa\kappa_\epsilon L(z)+\eta \Delta\}  \\
& \le 1 +  \frac{c}{\eta - c} \exp\{\kappa\kappa_\epsilon c_2 n +\eta \Delta\}. 
\end{align*}
The result follows. \qed
\\

{\bf Proof of Theorem \ref{SPLDP}} We first show that
\bes
\limsup_{n\to \infty} \frac{1}{n} \log p_n(z_n) \le -\inf_{x\in B} I(x).
\ees
From Proposition \ref{hitting}, it follows that 
\bes
\lim_{t\to\infty} \limsup_{n\to\infty} \frac{1}{n} \log \PP_{z_n}(\tau_n > t) \le \lim_{t\to\infty} \limsup_{n\to\infty} \frac{1}{n} \log \frac{\EE_{z_n}(e^{cn\tau_n})}{e^{cnt}} = -\infty,
\ees
where $\PP_{z_n}$ and $\EE_{z_n}$ are the probability and expectation conditioning on $Z_n(0)=z_n.$ Therefore we can take $T>0$ such that 
\bes
\limsup_{n\to\infty} \frac{1}{n} \log \PP_{z_n}(\tau_n > T)  < - \inf_{\phi\in F} I_{T,0}(\phi).
\ees
The following proof is similar to that of Lemma 5.7.21 in \cite{dembo2010}. Define the following closed set of $\mathcal{C}([0,T]; \RR^d_+)$:
\bes
F= \{\phi \in \mathcal{C}([0,T]; \RR^d_+): \psi \in \mathcal{C}([0,T]; \RR^d),  \phi(t) = \Gamma(\psi)(t)\in B \ \mbox{for some $t\in [0,T]$}\}.
\ees
Then 
\bes\ba
p_n(z_n) &= \PP_{z_n}(Z_n(\tau_n) \in B, \tau_n > T) + \PP_{z_n}(Z_n(\tau_n) \in B, \tau_n \le T) \\
&  \le \PP_{z_n}(\tau_n > T) + \PP_{z_n}(Z_n\in F).
\ea\ees
Also note that for $T>0$ we have
	$$
	\inf_{x\in B}I(x)\leq\inf_{\phi\in F}I_{T,0}(\phi)<\infty.
	$$
According to Lemma 1.2.16 in \cite{dembo2010}, 
\bes\ba
\limsup_{n\to \infty} \frac{1}{n} \log p_n(z_n) & \le \limsup_{n\to \infty} \frac{1}{n} \log (\PP_{z_n}(\tau_n > T) + \PP_{z_n}(Z_n\in F))\\
& \le \max \left[  \limsup_{n\to \infty} \frac{1}{n} \log \PP_{z_n}(\tau_n > T) ,  \limsup_{n\to \infty} \frac{1}{n}\log \PP_{z_n}(Z_n\in F) \right] \\
& \le - \inf_{\phi\in F} I_{T,0}(\phi) \le -\inf_{x\in B}I(x).
\ea\ees

Now we show the lower bound 
\bes
\liminf_{n\to\infty} \frac{1}{n} \log p_n(z_n) \ge -\inf_{\tilde z\in B} I(\tilde z) .
\ees
We next note that 
\bes
\inf_{\tilde z\in B} I(\tilde z) \ge \inf_{\psi\in\mathcal{AC}_0([0,\infty), \RR^d), \tau_{\psi} <\infty} \frac{1}{2} \int_0^{\tau_{\psi}} |\dot{\psi}(s)-\theta|^2_D ds,
\ees
where $\tau_{\psi} = \inf\{t\ge 0: \Gamma(\psi)(t) \in B\}.$  
Thus it suffices to establish that
\bes
\liminf_{n\to\infty} \frac{1}{n} \log p_n(z_n) \ge - \frac{1}{2} \int_0^{\tau_{\psi}} |\dot{\psi}(s)-\theta|^2_D ds
\ees
for any $\psi\in\mathcal{AC}([0,\infty), \RR^d)$ such that $\psi(0)=0$ and $\tau_{\psi} <\infty.$ Noting that the possible optimal solutions of the VP in \eqref{vp} never return to $0$ before reaching $B$, we only need to consider $\psi$ such that $\sum_{k=1}^d\Gamma(\psi)_k(s) > 0$ for all $s\in(0,\tau_\psi].$ 
Fix such a $\psi$. Let $\varphi^{\tilde z}(t) = \tilde z + \psi(t)$, with $\sum_{k=1}^d\tilde z_k>0.$
It is clear that 
\bes
\int_0^{\tau_\psi} |\dot{\varphi}^{\tilde z}(s)-\theta|_D^2 ds = \int_0^{\tau_\psi} |\dot{\psi}(s)-\theta|_D^2 ds.
\ees
Fix $\kappa > 0$, and let $B_\kappa = \{\tilde z \in \RR^d_+: \sum_{k=1}^d \tilde z_k = \kappa\}.$
For $\kappa, \delta > 0,$ define
\bes
G_\kappa(\varphi; \delta) = \cup_{\tilde z\in B_\kappa}\left\{\eta\in \mathcal{C}([0, \tau_\psi]; \RR^d_+): -\delta < (1-\delta)\sum_{k=1}^d \eta_k(t) - 
\sum_{k=1}^d \Gamma(\varphi^{\tilde z})_k(t) < \delta \ \mbox{for all $0\le t\le \tau_\psi$}\right\}.
\ees
For $\kappa\in (0,1)$ there exists $\delta_0>0$ and $N\in \NN$ such that, when $0<\delta\le \delta_0$, $$\sum_{i=1}^d \eta_i(t)> \left[\inf_{0\le t\le \tau_\psi} \sum_{i=1}^d \Gamma(\varphi^{\tilde z})_i(t) -\delta\right]/(1-\delta) > \epsilon/N,$$ for all $t\in[0,\tau_\psi]$, and $$\sum_{i=1}^d \eta_i(\tau_\psi) >\left[\sum_{i=1}^d \Gamma(\varphi^{\tilde z})_i(\tau_\psi) - \delta\right]/(1-\delta) \ge 1.$$

Assume that $n$ is large enough such that  $\kappa>\sum_kz_{n,k}>\epsilon/n$. From the Markov property of $Z_n$, we have that
\bes
p_n(z_n)  \ge \PP_{z_n}(Z_n \ \mbox{hits $B_\kappa$ before hitting $A_n$}) \PP\left(Z_n \ \mbox{hits $B$ before hitting $A_n$}|Z_n(0) \in B_\kappa\right).
\ees
Using the minimality property of the one-dimensional Skorokhod map, we see that $\sum_{k=1}^d Z_{n,k}(t) \ge \Gamma_1(\sum_{k=1}^d z_{n,k}+\sum_{k=1}^d \tilde X_{n,k})(t)$ for all $t\ge 0$, where $\tilde X_n(t) = \theta t + X(nt)/n$ and $\Gamma_1$ is the one-dimensional Skorokhod map. Thus we have 
\begin{align*}
&\PP_{z_n}(Z_n \ \mbox{hits $B_\kappa$ before hitting $A_n$}) \\
& \ge \PP_{\sum_{k=1}^d z_{n,k}}\left(\Gamma_1\left(\sum_{k=1}^d z_{n,k}+\sum_{k=1}^d \tilde X_k\right) \ \mbox{reaches $\kappa $ before reaching $\epsilon/n$}\right) \\
& =\PP_{\sum_{k=1}^d z_{n,k}}\left(\sum_{k=1}^d z_{n,k}+\sum_{k=1}^d \tilde X_k\ \mbox{reaches $\kappa$ before reaching $\epsilon/n$}\right)\\
& = \frac{1- e^{-n\gamma b_n}}{e^{n\gamma a_n}-e^{-n\gamma b_n}},
\end{align*}
where $\gamma = \frac{2|\sum \theta|}{\mathbf{1}'\Sigma \mathbf{1}}, a_n = \kappa -  \sum_{k=1}^d z_{n,k},$ and $b_n = \sum_{k=1}^d z_{n,k} - \frac{\epsilon}{n}$. Thus
\bes\ba
& \liminf_{n\to\infty} \frac{1}{n} \log \PP_{z_n}(Z_n \ \mbox{hits $B_\kappa$ before hitting $A_n$}) \\ & \ge  \liminf_{n\to\infty} \frac{1}{n} \log  \frac{1- e^{-n\gamma b_n}}{e^{n\gamma a_n}-e^{-n\gamma b_n}} \\
&   \ge \liminf_{n\to\infty} \left(\frac{\log(1-e^{-n \gamma b_n})}{n} - \gamma a_n \right)\\
& = -\gamma \kappa.
\ea\ees
We next note that 
\bes\ba
& \liminf_{n\to\infty} \frac{1}{n} \log \PP\left(Z_n \ \mbox{hits $B$ before hitting $A_n$}|Z_n(0) \in B_\kappa\right) \\
& \ge \liminf_{n\to\infty} \frac{1}{n} \log \inf_{Z_n(0)\in B_\kappa}\PP\left(Z_n \ \mbox{hits $B$ before hitting $A_n$}\right)  \\
& \ge  \liminf_{n\to\infty} \frac{1}{n} \log \inf_{Z_n(0)\in B_\kappa}\PP\left(Z_n \in G_\kappa(\varphi;\delta)\right) \\
& \ge - \sup_{x\in B_\kappa} \inf_{\eta \in G_\kappa(\varphi;\delta)} I_{\tau_\psi, x}(\eta) \\
& \ge  -\sup_{x\in B_\kappa} \frac{1}{2} \int_0^{\tau_{\psi}} \left|\frac{1}{1-\delta}\dot{\varphi^x}(s)-\theta\right|_D^2 ds \\
& = -\frac{1}{2} \int_0^{\tau_{\psi}} \left|\frac{1}{1-\delta}\dot{\psi}(s)-\theta\right|_D^2 ds.
\ea\ees 
Letting $\delta\downarrow 0$, we conclude that 
\[
\liminf_{n\to\infty} \frac{1}{n} \log p_n(z_n) \ge -\gamma\kappa  -\frac{1}{2} \int_0^{\tau_{\psi}} \left|\dot{\psi}(s)-\theta\right|_D^2 ds.
\]
Finally, letting $\kappa\downarrow 0$, the result follows. \qed
%{\color{red} Xin: Doesn't $\kappa$ need to go to 0 as well?} \col{Yes, corrected.} 

\subsection{Proofs for Section \ref{sec:vp}}

{\bf Proof of Lemma \ref{ss1}.} % We first note that for any $w, v\in S\backslash (A \cup B)$, there exists $K\subset I$ such that $w$ and $v$ satisfy Condition \ref{wvcondition} associated with $K$. 
For any $w, v\in \mathbb{R}^d_+\backslash (A \cup B)$, suppose we have an optimal triple $(\psi,\eta,\phi)$ with $\phi(0)=w, \phi(T)=v$. Denote by $\mathcal{P}(I)$ the power set of $L$. When $\phi$ traverses a face $F_K$ from time $t_1$ to $t_2$ for some $K\in \mathcal{P}(I)$, the segment of $(\psi,\eta,\phi)$ must be locally optimal. Otherwise, we can find a path which has lower cost. For some $K\in \mathcal{P}(I)$, we define the set $G_K=\{ t\in[0, T]: \phi_i(t)=0, \phi_j(t)>0, \forall i\in K, \forall j\in L\backslash K\}^{o}$, where for a set $E$, $E^o$ denotes its interior. We notice that $G_K$ is a open set, and can be represented as the disjoint union of at most countably many open intervals. Denote the set of these open intervals by $\mathbb{O}_K$. Notice that the optimal path traverses the face $F_K$ from time $t_1$ to $t_2$ if $(t_1,t_2)\in \mathbb{O}_K$. Observe that $\cup_{K\in \mathcal{P}(I)}\mathbb{O}_K$ contains at most countably many disjoint open intervals, and
\begin{align}
\mathcal{I}(w,v)&=\sum_{K\in \mathcal{P}(I)} \sum_{O\in \mathbb{O}_K} \mathcal{I}^*_{K}(\phi(\inf(O)),\phi(\sup(O)))\\
&\ge \sum_{K\in \mathcal{P}(I)} \sum_{O\in \mathbb{O}_K}\subsoln(\phi(\inf(O)))-\subsoln(\phi(\sup(O)))\\
&=\subsoln(w)-\subsoln(v),
\end{align}
Both of the equalities can be justified by the argument that there are at most countably many singleton points in the set $[0,T]\setminus \{\cup_{K\in \mathcal{P}(I)}G_K\}$, which contributes $0$ to the cost.
\qed
\\

{\bf Proof of Theorem \ref{constrain_vp}.} Theorem \ref{constrain_vp} summarizes some major results of \cite{el2012}, including Theorem 3, Theorem 4, Lemma 5, Proposition 1, and Proposition 2. Despite a small difference in problem settings, their proofs are still valid for our setting. To be more specific, given $K\subseteq L$, and two points $v,w \in \mathbb{R}^d_+$ which satisfy Condition \ref{wvcondition} with respect to $K$, Theorems 3 and 4 of \cite{el2012} assert that there is a unique $J_*\subseteq K$ such that  $\langle A^{J_*}R_j, \alpha^{J_*}(w,v)(v-w)-\theta \rangle_{D} \le 0, \  j\in K\backslash J_*$, and when $J_*\neq \emptyset$, each component of ${\lambda}^{J_*}(w,v)$ is positive. Proposition 1 of \cite{el2012} shows that $\mathcal{I}^*_K(w,v)=| A^{J_*}\theta |_{D}  | A^{J_*}(v-w) |_{D}-\langle A^{J_*}\theta , A^{J_*}(v-w)\rangle_{D}$. Lemma 5 and Proposition 2 of \cite{el2012} give the explicit form of an optimal path.       \qed
%{\color{red} Zicheng: Can you please add more details here?} 
\\

{\bf Proof of Lemma \ref{Kempty}.} See the proof of Lemma 1 in \cite{el2012}. \qed
\\

%{\color{red}I don't think this Lemma is necessary}
%{\bf Proof of Lemma \ref{Keq1}.} Assume $K=1$, $\gamma \langle R_1,v-w\rangle > \langle R_1,\theta \rangle$, and $J=\emptyset$. From Theorem \ref{constrain_vp}, 
%\begin{align*}
%\langle A^{\emptyset}R_1,\alpha(w,v)(v-w)-\theta\rangle \le 0.
%\end{align*}
%Since $A^{\emptyset}=\mathbb{I}_{d\times d}$, we have $\alpha(w,v)=\gamma$. Hence
%\begin{align*}
%\langle R_1, \gamma(v-w)-\theta) \rangle \le 0.
%\end{align*}
%which is contradictory to our assumption. To prove the reverse, we assume $K=\{1\}$, $J=\{1\}$, and $\gamma \langle %R_1,v-w\rangle \le \langle R_1,\theta \rangle$. However, this yields that 
%\begin{align*}
%	\langle A^{\emptyset}R_1,\alpha(w,v)(v-w)-\theta\rangle \le 0, 
%\end{align*}
%which implies that $J=\emptyset$ also satisfies Theorem \ref{constrain_vp}. From the uniqueness of such $J$, we reach a %contradiction to our assumption.  \qed
%\\

{\bf Proof of Lemma \ref{neg}.}
(i) If $\mathcal{I}^*_K(w,v)=0$, then by Theorem \ref{constrain_vp}, there exists an optimal path such that $\dot{\psi}(t)=\theta < 0$ for all $t\in (0, \tau_\phi(v))$, where $\tau_\phi(v) = \inf\{t\ge 0: \phi(t) = v\}$ given that $\phi(0) = w$. Hence
%{\color{red} Isn't this $\tau_\phi(v)$?} Hence 
\begin{equation}\label{diff_eqn}
\begin{aligned}
v-w & = \phi(\tau_\phi(v)) - \phi(0) = w + \theta\tau_\phi(v) +R\eta(\tau_\phi(v)) - w \\
& = \theta\tau_\phi(v) +R\eta(\tau_\phi(v)),
\end{aligned}
\end{equation}
where $\eta_i(\tau_\phi(v))>0,\forall i\in {K}$, and $\eta_i(\tau_\phi(v))=0,\forall i\in L\backslash {K}$ since $\phi(t)\in F_{K}$ for $t\in (0, \tau_\phi(v))$. From Theorem \ref{constrain_vp},
%{\color{red} Why cite Proposition 2 if you already cited in statement of Theorem \ref{constrain_vp}?}and Proposition 2 of \cite{el2012}
we have
\begin{equation}\label{diff_eqn2}
b_K(w,v)=\alpha^{J_*}(w,v)(v-w)-R_{LJ_*}\alpha^{J_*}(w,v)B^{J_*}(v-w)+R_{LJ_*}B^{J_*}\theta=\theta,
\end{equation}
and
\begin{equation}\label{diff_eqn3}
\tau_\phi(v)=\frac{1}{\alpha^{J_*}(w,v)}.
\end{equation}
Multiply both sides of the second equality in \eqref{diff_eqn2} by $\tau_\phi(v)$, and we can get
\begin{equation*}
v-w=R_{LJ_*}B^{J_*}(v-w)-R_{LJ_*}B^{J_*}\theta\tau_\phi(v)+\theta\tau_\phi(v).
\end{equation*}
Comparing the above equation with \eqref{diff_eqn}, we obtain that
\begin{equation*}
R\eta(\tau_\phi(v))=R_{IJ_*}B^{J_*}(v-w-\theta\tau_\phi(v)).
\end{equation*}
Since $R$ is invertible, all its columns are linearly independent. If $J_*\neq K$,  for $i\in K\setminus J_*$, $\eta_i(\tau_\phi(v)) =0$, which leads to a contradiction. Thus, we can conclude that $J=K$.

%Plug \ref{diff_eqn2} and \ref{diff_eqn3} into \ref{diff_eqn}, we can observe that $J=K$. {\color{red} I think this %argument needs more detail.}

We next note that $v_i - w_i < 0$ for $i\in L\backslash K$ from \eqref{diff_eqn} and the fact that $\theta < 0$; which yields that $v_i - w_i\le 0,$ for all $i \in L$. Noting that both $w$ and $v$ are in $\RR^d_+$, $\norm(v)\le \norm (w)$ follows from the monotonicity of $\norm(\cdot)$. 

(ii) 
Again from Theorem \ref{constrain_vp}, we can find two optimal paths such that $\dot{\psi}_1(t)=\theta, t\in(0, \tau(v))$ and  $\dot{\psi}_2(t)=\theta, t\in(0, \tau(w))$, where $\tau(v) = \inf\{t\ge 0: \phi_1(t) = v\}$ and $\tau(w) = \inf\{t\ge 0: \phi_2(t) = w\}$ given that $\phi_1(0)=\phi_2(0) = 0$. For simplicity, we ignore the symbols of particular optimal path in the following proof. Similar to \eqref{diff_eqn}, we have 
\begin{align*}
v & = \theta\tau(v) +R\eta(\tau(v)), \  \ w  = \theta \tau(w) + R\eta(\tau(w)),
\end{align*}
where $\eta_i(\tau(v))>0, \eta_i(\tau(w)) > 0,\forall i\in {K}$, and $\eta_i(\tau(v))=\eta_i(\tau(w)) = 0,\forall i\in L\backslash {K}$. Let $\phi_K = (\phi_i; i\in K), \eta_K = (\eta_i; i\in K),$ and $\theta_K = (\theta_i; i\in K)$. We first consider $\mathcal{I}^*_K(v)$, noting that $\eta_i(t) = 0$ for $i\in L\backslash K$ and $t\in (0, \tau(v))$, we have for $t\in (0, \tau(v)]$,
\[
0 = \phi_K(t) =   \theta t + R_{KK}\eta_K(t), 
\]
which gives that $\eta_K (t) = - R^{-1}_{KK} \theta t.$ Similarly, for $\mathcal{I}^*_K(w)$, we have $\eta_K (t) = - R^{-1}_{KK} \theta t$ for $t\in (0, \tau(w)].$ It now follows that for $i\in K,$
\begin{align*}
v_i & = \tau(v)(\theta_i  - [R^{-1}_{KK}\theta]_i),\ \ w_i =  \tau(w)(\theta_i - [R^{-1}_{KK}\theta]_i),
\end{align*}
and for $i\in L\backslash K,$
\begin{align*}
v_i & = \tau(v)\theta_i,\ \ w_i  =  \tau(w)\theta_i.
\end{align*}
The result follows. 
\qed 
\\

%{\bf Proof of Lemma \ref{lem:neg}}
%By Lemma \ref{neg}, if $\mathcal{I}^*_K(w,v)=0$, then $u_i\le 0, \forall i$, where $u=v-w$. Therefore $W(v)\le W(w)$ follows from the monotonicity of $W(\cdot)$. \qed
%\\
%
%{\bf Proof of Lemma \ref{lem:linear}.} Use the similar argument in proving \ref{neg}, $v=k_1(\theta+RY)$, $w=k_2(\theta+RY)$ for some negative numbers $k_1$ and $k_2$. Therefore, $v=cw$ for some constant $c$. \qed
%\\

 {\bf Proof of Proposition \ref{prop:pos}.} It suffices to show that for any $K\subset L$, there exists an $\epsilon>0$ such that for any $v\in \mathbb{D}_K$, $\mathcal{I}^*_K(v) > \epsilon$. From the extreme value theorem, $\mathcal{I}^*_K(v)$ attains the minimum in $\mathbb{D}_K$. It is also known that $\mathcal{I}^*_K(v)\neq 0$ for any $v\in \mathbb{D}_K$. To see this, if there exists a $v\in\mathbb{D}_K$ such that $\mathcal{I}^*_K(v) =0$, then from Lemma \ref{neg} (i), $v_j<0$ for all $j\in L\setminus K$, which contradicts the property (3) in the definition of $\mathbb{D}_K$. Hence there exists some positive number $\epsilon$ such that  $\mathcal{I}^*_K(v) > \epsilon$ for any $v\in \mathbb{D}_K$.\qed
\\

{\bf Proof of Theorem \ref{th:sub-solution}.} From Lemma \ref{ss1}, it suffices to show that $\subsoln(x)\le 0$ for all $x\in B$, $\subsoln$ is continuous, and that for all $w,v  \in \mathbb{R}^d_+\backslash B$, which satisfy Condition \ref{wvcondition} for some $K\subset L$, we have 
$$
\subsoln(w)-\subsoln(v)\leq\mathcal{I}^*_K(w,v).
$$
Since the first two conditions are obviously satisfied we only need establish the third condition.
First assume $v_i - w_i <0, \forall i\in L\backslash K$, then
\begin{align*}
\subsoln(w)-\subsoln(v)& =\frac{{\norm}(v)-{\norm}(w)}{\mathbf{r}} \le 0 \le \mathcal{I}^*_{K}(w,v).
\end{align*} 
Now if there exists at least one $j\in L\backslash K$ such that $v_j - w_j \ge 0$, we have $(v-w)/|v-w| \in \mathbb{D}_K$, and 
\begin{align*}
\subsoln(w)-\subsoln(v)&=\frac{{\norm}(v)-{\norm}(w)}{\mathbf{r}} \\
& \le \frac{|v-w|\norm((v-w)/|v-w|)}{\mathbf{r}}\\
& \le |v-w|\mathcal{I}^*_{K}((v-w)/|v-w|)\\
& = \mathcal{I}^*_{K}(v-w),
\end{align*} 
where the first inequality follows since $T$ is a norm, the second inequality follows from the definition of $\mathbf{r}$ and the final equality follows from Theorem \ref{constrain_vp}.\qed

\subsection{Proofs for Section \ref{sec:splitting}}

\subsubsection{Proofs of Propositions \ref{exp_mome} and \ref{error}}

{\bf Proof of Proposition \ref{exp_mome}.} Let $W$ be a standard Brownian motion, and for a fixed $\Delta>0$, define its discretization $W^\Delta(t) = W((k-1)\Delta)$ for $t\in[(k-1)\Delta, k\Delta)$ and $k\in\NN.$ Let $T>0$, and consider $\mathbf{W}_T=\{W(t):0\leq t\leq T\}$ and $\mathbf{W}_T^\Delta=\{W^\Delta(t): 0\leq t\leq T\}$. 
%Define $\mathbf{W}_T=\{W(t):0\leq t\leq T\}$ to be a standard Brownian motion. For a fixed $\Delta>0$ we will consider the discretization $\mathbf{W}_T^\Delta=\{W^\Delta(t): 0\leq t\leq T\}$ defined by $W^\Delta(t)=W((k-1)\Delta)$ for $t\in [(k-1)\Delta, k\Delta)$ and $k\in\NN$. 
We are interested in the following discretization error.
\begin{align*}
d(\mathbf{W}_T,\mathbf{W}_T^\Delta)&=\max_{0\leq t\leq T}|W(t)-W^\Delta(t)|\\
&\le
\max_{0\leq k\leq \left \lfloor{T/\Delta}\right \rfloor }\max_{t\in [k\Delta, (k+1)\Delta)}|W(t)-W(k\Delta)|.
\end{align*}
Define $V(\Delta)=\max_{0\leq t\leq \Delta}|W(t)|$. By the Markov property of Brownian motion, we have that
$$
d(\mathbf{W}_T,\mathbf{W}_T^\Delta)\le\max_{1\leq k\leq \left \lfloor{T/\Delta}\right \rfloor }V_k(\Delta),
$$
where $\{V_k(\Delta): k\in \{1,\ldots, \left \lfloor{T/\Delta}\right \rfloor\}$ is an independent sequence of copies of $V(\Delta)$. From \cite{karatzas2012brownian}, an upper bound on the tail probabilities of $V(\Delta)$ can be derived. In particular, for $z>0$,
$$
\PP(V(\Delta)>z)\leq \frac{4}{z}\sqrt{\frac{\Delta}{2\pi}}e^{-z^2/(2\Delta)}.
$$
%For ease of notation define $X_T(\Delta)=d(\mathbf{W}_T,\mathbf{W}_T^\Delta)$ and 
It follows that for $z>0$
\begin{align}
\label{eq:UpperBoundTailProbV}
\PP(d(\mathbf{W}_T,\mathbf{W}_T^\Delta)>z)& = 1-\PP(V(\Delta)\leq z)^{T(\Delta)}\leq 1-\left(1-\frac{4}{z}\sqrt{\frac{\Delta}{2\pi}}e^{-z^2/(2\Delta)}\right)^{T(\Delta)}
\nonumber\\
&=
1-\exp\left[ T(\Delta)\log\left(1-\frac{4}{z}\sqrt{\frac{\Delta}{2\pi}}e^{-z^2/(2\Delta)}\right)\right],
\end{align}
where $T(\Delta)=\left \lfloor{T/\Delta}\right \rfloor$.
Next we use the inequality that $\log (z)\geq -(1-z)/z$ for $z\in (0,1)$ to conclude that
$$
\PP(d(\mathbf{W}_T,\mathbf{W}_T^\Delta)>z) \leq 1 - \exp\left(-T(\Delta)\left(\frac{\frac{4}{z}\sqrt{\frac{\Delta}{2\pi}}e^{-z^2/(2\Delta)}}{1-\frac{4}{z}\sqrt{\frac{\Delta}{2\pi}}e^{-z^2/(2\Delta)}}\right)\right).
$$
Considering the exponent in the above equation, we have 
\begin{align*}
-T(\Delta)\left(\frac{\frac{4}{z}\sqrt{\frac{\Delta}{2\pi}}e^{-z^2/(2\Delta)}}{1-\frac{4}{z}\sqrt{\frac{\Delta}{2\pi}}e^{-z^2/(2\Delta)}}\right) & = - \left(\frac{T}{\Delta}+O(1)\right) \frac{4\sqrt{\Delta}}{\sqrt{2\pi}z e^{z^2/(2\Delta)}- 4\sqrt{\Delta}} \\
& = - (1+O(\Delta))\frac{4T}{\sqrt{2\pi\Delta}z e^{z^2/(2\Delta)}- 4\Delta},
\end{align*}
and thus 
\begin{align}
\label{eq:XT_upperbound}
\PP(d(\mathbf{W}_T,\mathbf{W}_T^\Delta)>z)\leq 1-\exp\left( - \frac{(1+O(\Delta))4T}{\sqrt{2\pi\Delta}z e^{z^2/(2\Delta)}- 4\Delta}\right)\leq \frac{(1+O(\Delta))4T}{\sqrt{2\pi\Delta}z e^{z^2/(2\Delta)}- 4\Delta}.
\end{align}

Now consider a scaled version of $\mathbf{W}_T$ defined by $W_n(t)=\frac{1}{n}W(nt), t\in [0, T]$, which is denoted  by $\mathbf{W}_{n,T}=\{W_n(t): 0\leq t\leq T\}$. For $\Delta>0$ we consider the discretization $\mathbf{W}^\Delta_{n,T}=\{W_n^\Delta(t): 0\leq t\leq T\}$ defined by $W_n^\Delta(t)=W_n((k-1)\Delta)$ for $t\in [(k-1)\Delta, k\Delta).$ Finally, we observe that
$$
d(\mathbf{W}_{n,T},\mathbf{W}_{n,T}^\Delta) =  \frac{1}{n}d(\mathbf{W}_{nT},\mathbf{W}_{nT}^{n\Delta}).
%X_T^n(\Delta)=d(\mathbf{W}_{n,T},\mathbf{W}_{n,T}^\Delta),
$$
%and observe that $X_T^n(\Delta)=\frac{1}{n}X_{nT}(n\Delta)$.
Therefore we can use \eqref{eq:XT_upperbound} to obain 
\begin{align}\label{eq:scaled_XT_upperbound}
\PP(d(\mathbf{W}_{n,T},\mathbf{W}_{n,T}^\Delta)>z)=\PP(d(\mathbf{W}_{nT},\mathbf{W}_{nT}^{n\Delta})>nz)\leq  \frac{(1+O(n\Delta))4nT}{\sqrt{2\pi n\Delta} \ nz e^{nz^2/(2\Delta)}- 4n\Delta}.
\end{align}

Under Assumption \ref{m_matrix}, the Skorokhod map is Lipschitz continuous, and there exists a constant $C$ which only depends on $R$, such that for $t\ge 0$, 
\begin{align}
\sup_{0\le s \le t}|Z(s) - Z^\Delta(s)| \le C \sup_{0\le s \le t}|\tilde X(s) - \tilde X^\Delta(s)|.
\end{align}
For $z>0$, we have 
\begin{align*}
\PP\left(\sup_{0\le s \le t}|Z(s) - Z^\Delta(s)| > z\right) & \le \PP\left(\sup_{0\le s \le t}|\tilde X(s) - \tilde X^\Delta(s)| > z/C\right) \\
& \le \PP\left(\sup_{0\le s \le t}|X(s) - X^\Delta(s)| > z/C - \Delta|\theta|\right) \\
& \le \PP\left(\sup_{0\le s \le t}|W(s) - W^\Delta(s)| > \frac{z/C - \Delta|\theta|}{d|\Sigma^{1/2}|}\right) 
\end{align*}
Denote $\tilde z_\Delta =  (z/C - \Delta|\theta|)/(d|\Sigma^{1/2}|).$ From \eqref{eq:XT_upperbound}, we have 
\begin{align}
\PP\left(\sup_{0\le s \le T}|Z(s) - Z^{\Delta}(s)| > z\right) & \le \PP\left(\sup_{0\le s \le T}|W(s) - W^\Delta(s)| > \tilde z_\Delta\right) \label{eq:prop5.1} \\
& \le  \frac{(1+O(\Delta))4T}{\sqrt{2\pi\Delta}\tilde z_\Delta e^{\tilde z_\Delta^2/(2\Delta)}- 4\Delta}.\nonumber
\end{align}

Now consider the large derivation scaled processes $Z_n(t) = Z(nt)/n, X_n(t) = X(nt)/n$, and $Y_n(t) = Y(nt)/n$. From the homogenity property of the Skorokhod map for linear scalings, $(Z_n, Y_n)$ is the solution of the Skorokhod problem for $\{\tilde X_n(t) = X_n(t) + \theta t; t\ge 0\}$. It follows from \eqref{eq:scaled_XT_upperbound} that 
\begin{align*}
\PP\left(\sup_{0\le s \le T}|Z_n(s) - Z_n^\Delta(s)| > z\right) & \le \PP\left(\sup_{0\le s \le T}|W_n(s) - W_n^\Delta(s)| > \tilde z_\Delta\right) \\
& \le \frac{(1+O(n\Delta))4nT}{\sqrt{2\pi n\Delta} \ n\tilde z_\Delta e^{n\tilde z_\Delta^2/(2\Delta)}- 4n\Delta}.
\end{align*}
Hence 
\begin{align}
& \frac{1}{n}\log \PP\left(\sup_{0\le s \le T}|Z_n(s) - Z_n^\Delta(s)| > z\right) \label{exp_equiv}\\
& \leq  \frac{1}{n}\log\left( \frac{(1+O(n\Delta))4nT}{\sqrt{2\pi n\Delta} \ n\tilde z_\Delta e^{n\tilde z_\Delta^2/(2\Delta)}- 4n\Delta}\right) \nonumber\\
& = \frac{\log((1+O(n\Delta))4nT)}{n} - \frac{1}{n}\log\left(\sqrt{2\pi n\Delta} \ n\tilde z_\Delta e^{n\tilde z_\Delta^2/(2\Delta)}- 4n\Delta\right) \nonumber \\
& =  \frac{\log((1+O(n\Delta))4nT)}{n}  -  \frac{1}{n} \left[\log(e^{n\tilde z_\Delta^2/(2\Delta)}) + \log(n) + \log(\sqrt{2\pi n \Delta} \tilde z_\Delta - 4 \Delta e^{-n\tilde z_\Delta^2/(2\Delta)})\right] \nonumber\\
& =  \frac{\log((1+O(n\Delta))4nT)}{n} -\frac{\tilde z_\Delta^2}{2\Delta}  - \frac{1}{n} \left[\log(n) + \log(\sqrt{2\pi n \Delta} \tilde z_\Delta - 4 \Delta e^{-n\tilde z_\Delta^2/(2\Delta)})\right]. \nonumber
\end{align}
When $\Delta = \Delta(n) = o(1)$, we have 
\[
\limsup_{n\to\infty} \frac{1}{n}\log \PP\left(\sup_{0\le s \le T}|Z_n(s) - Z_n^\Delta(s)| > z\right) = -\infty.
\]
\qed
\\

To show Proposition \ref{error}, we introduce the following notation. Let $\eta \in (0, 1/2)$, and take $\Delta\equiv\Delta(n)$ small enough such that $\epsilon/n-{\Delta}^{\eta} > 0$ and $\epsilon/n+{\Delta}^{\eta} < z_n$.

\begin{lemma}\label{exp_mome2}
	If $\Delta(n)\rightarrow 0$ as $n\rightarrow \infty$, then for $\eta \in (0, 1/2)$,
	\begin{equation}
	\limsup_{n\to\infty}\frac{1}{n}\log \PP \left( \sup_{0\le s \le \tau_n} |Z_n(s)-Z_n^{\Delta}(s)|>{\Delta}^{\eta}\right) = -\infty.
	\end{equation}
\end{lemma}
{\bf Proof.} From \eqref{eq:prop5.1} in the proof of Proposition \ref{exp_mome}, letting $\xi = \frac{{{\Delta}^{\eta}}/C - \Delta|\theta|}{d|\Sigma^{1/2}|}$, we have that for $T\ge 0,$
\begin{align*}
& \quad \quad \frac{1}{n}\log \PP \left( \sup_{0\le s \le T} |Z_n(s)-Z_n^{\Delta}(s)|>{\Delta}^{\eta}\right)\\
& \le  \frac{\log(4nT)}{n} -\frac{ {\xi}^2}{2\Delta}  - \frac{1}{n} \left[\log(n) + \log(\sqrt{2\pi n \Delta} \ {\xi}- 4 \Delta e^{-n{\xi}^2/(2\Delta)})\right].
%& = \frac{\log(4nT)}{n} -\frac{(\frac{{\xi}/C - \Delta}{|\Sigma^{1/2}|})^2}{2\Delta} - \frac{1}{n} \left[\log(n) + \log(\sqrt{2\pi n \Delta}\frac{{\xi}/C - \Delta}{|\Sigma^{1/2}|} - 4 \Delta e^{-n(\frac{{\xi}/C - \Delta}{|\Sigma^{1/2}|})^2/(2\Delta)})\right]
\end{align*}
Next noting that 
\[
\lim_{n\to\infty}\frac{\xi^2}{2\Delta}=\infty,
\]
%When $\xi$ is substituted by ${\Delta}^{(1/2-\eta)}$, $\lim_{n\to\infty}\frac{(\frac{{{\Delta}^{(1/2-\eta)}}/C - \Delta}{|\Sigma^{1/2}|})^2}{2\Delta}=\infty$. 
it follows that
\begin{equation}\label{finite_time_est}
\limsup_{n\to\infty}\frac{1}{n}\log \PP \left( \sup_{0\le s \le T} |Z_n(s)-Z_n^{\Delta}(s)|>{\Delta}^{\eta}\right)=-\infty.
\end{equation}
We next observe that for $T>0$
	\begin{align*}
	& \PP \left( \sup_{0\le s \le \tau_n} |Z_n(s)-Z_n^{\Delta}(s)|>{\Delta}^{\eta}\right)\\
	&\leq \PP\left(\sup_{0\le s \le T} |Z_n(s)-Z_n^{\Delta}(s)|>{\Delta}^{\eta}\right)+ \PP\left( \tau_n> T \right).
	\end{align*}
	and therefore by Lemma 1.2.16 in \cite{dembo2010}
	\begin{align*}
	&\limsup_{n\to\infty}\frac{1}{n}\log\PP\left( \sup_{0\le s \le \tau_n} |Z_n(s)-Z_n^{\Delta}(s)|>{\Delta}^{\eta}\right)\\
	&\leq
	\limsup_{n\to\infty}\frac{1}{n}\log\PP\left(\tau_n>T\right).
	\end{align*}
	The result then follows from Proposition \ref{hitting} and sending $T\to\infty$. \qed

%\begin{lemma}\label{sandwich}
%	If $\Delta(n)\rightarrow 0$ as $n\rightarrow \infty$, then
%	\begin{equation*}
%	\limsup_{n\rightarrow \infty} \frac{p_n^{\Delta}(z_n)}{p_n^-(z_n)} \le 1+\limsup_{n\rightarrow \infty} \frac{\PP\left(\sup_{0\le s \le \tau_n^-}|Z_n(s) - Z_n^\Delta(s)| > {\Delta}^{\eta}\right)}{p_n^-(z_n)}.
%	\end{equation*}
%\end{lemma}
%{\bf Proof.} Let $p_n^-(z_n)$ denote the probability that the SRBM $Z_n$ reaches $\{\tilde z \in \RR^d_+: \sum_{k=1}^d \tilde z_k \ge 1-\Delta^\eta\}$ before reaching $\{\tilde z \in \RR^d_+: \sum_{k=1}^d \tilde z_k \le \epsilon/n - \Delta^\eta\}$ given that $Z_n(0) = z_n$. We observe that 
%\begin{align*}
%p_n^{\Delta}(z_n)& = \PP\left(\tau_n^{\Delta}(B_n)<\tau_n^{\Delta}(A_n),\sup_{0\le s \le \tau_n^-}|Z_n(s) - Z_n^\Delta(s)| < {\Delta}^{\eta}\right)\\
%& \quad \quad +\PP\left(\tau_n^{\Delta}(B_n)<\tau_n^{\Delta}(A_n),\sup_{0\le s \le \tau_n^-}|Z_n(s) - Z_n^\Delta(s)| > {\Delta}^{\eta}\right)\\
%& \le p_n^{-}(z_n)+\PP\left(\sup_{0\le s \le \tau_n^-}|Z_n(s) - Z_n^\Delta(s)| > {\Delta}^{\eta}\right)
%\end{align*}
%Hence
%\begin{align*}
%\frac{p_n^{\Delta}(z_n)}{p_n^-(z_n)}-1\le \frac{\PP\left(\sup_{0\le s \le \tau_n^-}|Z_n(s) - Z_n^\Delta(s)| > {\Delta}^{\eta}\right)}{p_n^-(z_n)}.
%\end{align*} \qed\\

%{\color{red} Lemma \ref{sandwich} is pretty trivial. Perhaps we remove it? {\color{blue} Removed!}} \\

Define:
\begin{align*}
& A_n^-=\left\{z\in \RR^d_+: \sum_{k=1}^d z_k  \le \epsilon/n-{\Delta}^{\eta}\right\},\ \ \mbox{and} \ \ B_n^- = \left\{z\in \RR^d_+: \sum_{k=1}^d z_k \ge 1-{\Delta}^{\eta}\right\},
%& A_n^+=\{z\in \RR^d_+: \sum_{k=1}^d z_k  \le \epsilon/n+{\Delta}^{\eta}\},\ \ \mbox{and} \ \ B_n^+ = \{z\in \RR^d_+: \sum_{k=1}^d z_k \geq 1+{\Delta}^{\eta}\}.
\end{align*}
and 
\begin{align*}
& \tau_n^-(A_n^-) = \inf\{t\geq 0: Z_n(t)\in A_n^-\}, \  \tau_n^-(B_n^-) = \inf\{t\geq 0: Z_n(t)\in  B_n^-\},\\
&\tau_n^- = \inf\{t\geq 0: Z_n(t)\in A_n^-\cup B_n^-\}.
\end{align*}
At last, for $\tilde z\in \RR^d_+$, let $p_n^-(\tilde z) = \PP(\tau_n^-(B_n^-) < \tau_n^-(A_n^-)| Z_n(0) = \tilde z).$ 
\begin{lemma}\label{san2}
	For $\eta\in (0, 1/2)$, let $\Delta \equiv \Delta(n)$ satisfy $\lim\limits_{n\rightarrow \infty}n\Delta^{\eta}=0$. Then 
	\[\frac{\max_{\tilde z\in \initial_n}p_n^-(\tilde z)}{\max_{\tilde z\in \initial_n} p_n(\tilde z)}\le 1+O\left(n\Delta^{\eta}\right),\]
	where $\initial_n = \{\tilde z\in \RR^d_+: \sum_{i=1}^d \tilde z_i = \sum_{i=1}^d z_{n,i}\}.$
\end{lemma}
{\bf Proof.} Define the boundaries of $A^-_n, B^-_n$ and $A_n$ inside $\RR^d_+$ as follows:
\begin{align*}
& \cla_n^-=\left\{z\in \RR^d_+: \sum_{k=1}^d z_k  = \epsilon/n-{\Delta}^{\eta}\right\}, \clb_n^- = \left\{z\in \RR^d_+: \sum_{k=1}^d z_k = 1-{\Delta}^{\eta}\right\}, 
\end{align*}
and
\[
 \cla_n=\left\{z\in \RR^d_+: \sum_{k=1}^d z_k  = \epsilon/n\right\}.
\]
Using the Markov property,  
\begin{align*}
& p_n^-(z_n) - p_n(z_n) \\
& = \PP(\tau_n^-(B^-_n) < \tau_n^-(A_n^-), \tau_n(B) > \tau_n(A_n) | Z_n(0) =z_n) \\
& = \PP(\tau_n^-(B^-_n) < \tau_n(A_n) < \tau_n(B)  | Z_n(0) =z_n) \\
& \quad +  \PP(\tau_n(A_n) < \tau_n^-(B^-_n) < \tau_n^-(A_n^-) | Z_n(0) =z_n) \\
& \le \max_{\tilde z\in \clb^-_n}\PP(\tau_n(A_n) < \tau_n(B)  | Z_n(0) = \tilde z) \PP(\tau_n^-(B^-_n) < \tau_n(A_n) | Z_n(0)=z_n) \\
& \quad +  \max_{\tilde z\in \cla_n}\PP(\tau_n^-(B^-_n) < \tau_n^-(A_n^-) | Z_n(0) = \tilde z)  \PP(\tau_n(A_n) < \tau_n^-(B^-_n)| Z_n(0) =z_n). 
\end{align*}
We observe that $\{\tau^-_n(B_n^-) < \tau_n(A_n)|Z(0)=z_n\}\subset \{\tau_n^-(B^-_n) < \tau_n^-(A_n^-)|Z(0)=z_n\}$, and thus
\[
\PP(\tau_n^-(B^-_n) < \tau_n(A_n) | Z_n(0)=z_n) \le p^-_n(z_n).
\]
Next note that for $\tilde z\in \cla_n,$
\begin{align*}
& \max_{\tilde z \in \cla_n}\PP(\tau_n^-(B^-_n) < \tau_n^-(A_n^-) | Z_n(0) = \tilde z) \\
& \leq \max_{\tilde z \in \initial_n}p^-_n(\tilde z)\max_{\tilde z \in \cla_n}\PP(\tau_n(C_n) < \tau_n^-(A_n^-) | Z_n(0) = \tilde z),
\end{align*}
where $\tau_n(\initial_n) = \inf\{t\ge 0: Z_n(t)\in \initial_n\}$.
It follows that
\begin{align*}
p_n^-(z_n) - p_n(z_n) & \le p^-_n(z_n) \max_{\tilde z\in \clb^-_n}\PP(\tau_n(A_n) < \tau_n(B)  | Z_n(0) = \tilde z) \\
& \quad + \max_{\tilde z \in \initial_n}p^-_n(\tilde z)\max_{\tilde z \in \cla_n}\PP(\tau_n(C_n) < \tau_n^-(A_n^-) | Z_n(0) = \tilde z), 
\end{align*}
which yields that
\begin{align*}
& p_n^-(z_n) [1- \max_{\tilde z\in \clb^-_n}\PP(\tau_n(A_n) < \tau_n(B)  | Z_n(0) = \tilde z) ] \\
& \le p_n(z_n) +  \max_{\tilde z \in \initial_n}p^-_n(\tilde z)\max_{\tilde z \in \cla_n}\PP(\tau_n(C_n) < \tau_n^-(A_n^-) | Z_n(0) = \tilde z).
\end{align*}
Taking maximum over $z_n \in \initial_n$, we have 
\begin{align*}
& \max_{\tilde z\in \initial_n}p_n^-(\tilde z) [1- \max_{\tilde z\in \clb^-_n}\PP(\tau_n(A_n) < \tau_n(B)  | Z_n(0) = \tilde z) ] \\
& \le \max_{\tilde z \in \initial_n} p_n(\tilde z) +  \max_{\tilde z \in \initial_n}p^-_n(\tilde z)\max_{\tilde z \in \cla_n}\PP(\tau_n(C_n) < \tau_n^-(A_n^-) | Z_n(0) = \tilde z),
\end{align*}
and 
\begin{align}
& \frac{\max_{\tilde z\in \initial_n}p_n^-(\tilde z)}{\max_{\tilde z \in \initial_n} p_n(\tilde z)}  \le \nonumber\\
& \frac{1}{1- \max_{\tilde z\in \clb^-_n}\PP(\tau_n(A_n) < \tau_n(B)  | Z_n(0) = \tilde z)-\max_{\tilde z \in \cla_n}\PP(\tau_n(C_n) < \tau_n^-(A_n^-) | Z_n(0) = \tilde z)}.\label{ineq:7.4.4}
\end{align}
Similar to the proof of Theorem \ref{SPLDP}, using the minimality property of the one-dimensional Skorokhod map $\Gamma_1$, we have 
\begin{align}
& \max_{\tilde z\in \clb^-_n}\PP(\tau_n(A_n) <  \tau_n(B) |Z_n(0)=\tilde z)\nonumber \\
& \le  \PP\left(\Gamma_1\left(1- \Delta^\eta+\sum_{k=1}^d \tilde X_k\right)  \mbox{reaches $\epsilon/n$ before reaching $1$}\right)\\
& =  \PP\left(1 - \Delta^\eta +\sum_{k=1}^d \tilde X_k  \ \mbox{reaches $\epsilon/n$ before reaching $1$}\right)\nonumber\\
& = \frac{1- e^{-n\gamma \Delta^\eta}}{e^{n\gamma a_n}-e^{-n\gamma  \Delta^\eta}},\label{ineq:7.4.3}
\end{align}
where $\gamma = \frac{2|\sum \mu|}{\mathbf{1}'\Sigma \mathbf{1}}$ and $a_n = 1 - \Delta^\eta - \epsilon/n.$ 

We next use the Lyapunov function $L$ introduced in Theorem \ref{lyp} to derive an upper bound for $\max_{\tilde z \in \cla_n}\PP(\tau_n(C_n) < \tau_n^-(A_n^-) | Z_n(0) = \tilde z)$. Choose $\kappa_1 > 0$ such that when $\sum_{i=1}^d z_i \ge \kappa_1$, $tr(\ddot{L}(z)\Sigma) < c_0$. Let $\kappa_2 = 2\kappa_1/\epsilon$, and define $\check Z(t) = \kappa_2 Z(t)$ for $t\ge 0$. Then for $\tilde z \in \cla_n$, 
\begin{align*}
\p(\tilde z)&\equiv  \PP(\tau_n(C_n) < \tau_n^-(A_n^-) | Z_n(0) = \tilde z) \\
& = \PP(\check Z(\cdot) \ \mbox{first reaches $\kappa_2 n\initial_n$ before reaching $\kappa_2 n \cla^-_n$}| \check Z(0) = \kappa_2 n \tilde z).
\end{align*} 
%We next let $\underline{l} = \inf\{L(z): z\in\RR^d_+, \sum_{i=1}^dz_i =1\}$, and then $\underline l \in (0, \infty)$. Fix $\tilde z\in A_n$, and let $z_0 = \kappa_2 n \tilde z$. Define a new Lyapunov function
%\[
%L_{\tilde z}(\check z) = \frac{\kappa_2\epsilon \underline{l}   }{L(z_0)} L(\check z), \ \check z\in\RR^d_+.
%\]
From Theorem \ref{lyp} (v),  for $\check z\in \RR^d_+,$ $L(\check z) = (\sum_{i=1}^d \check z_i) L(\check z/\sum_{i=1}^d \check z_i)$. Let $\underline{l} = \inf\{L(z): z\in\RR^d_+, \sum_{i=1}^dz_i =1\}$. Then $\underline l \in (0, \infty)$, and 
\begin{align}\label{ineq:7.4.1}
\sum_{i=1}^d \check z_i \le \frac{L(\check z)}{\underline l}. % \le (\kappa_3\vee 1) L(\check z). 
\end{align}
We further define the following two sets
\begin{align*}
\tilde A & = \left\{\check z\in \RR^d_+: L(\check z) \le \underline{l}\kappa_2(\epsilon-n\Delta^\eta)\right\},\\
\tilde C & = \left\{\check z\in \RR^d_+: L(\check z) \ge \underline{l}\kappa_2\sum_{i=1}^d z_i \right\}.
\end{align*}
We will see that $\tilde A \subset \kappa_2 n A^-_n$ and $ \kappa_2 n C_n\subset \tilde C$, where $C_n = \{\check z \in \RR^d_+: \sum_{i=1}^d \check z_i \ge \sum_{i=1}^d z_{n,i}\}.$ Indeed, for $\check z \in \tilde A$, from \eqref{ineq:7.4.1}, 
\begin{align*}
\sum_{i=1}^d \check z_i \le  \frac{L(\check z)}{\underline{l}} \le \underline{l}\frac{\kappa_2(\epsilon-n\Delta^\eta)}{\underline{l}} = \kappa_2 n(\epsilon/n - \Delta^\eta), 
\end{align*}
which says $\check z \in \kappa_2 n A^-_n$. Now for $\check z \in \kappa_2 n \mathcal{C}_n$, that is $\sum_{i=1}^d \check z_i \ge \kappa_2 \sum_{i=1}^d z_i$, again from \eqref{ineq:7.4.1}, we have 
\begin{align*}
L(\check z) \ge \underline{l}\sum_{i=1}^d \check z_i \ge \underline{l} \kappa_2\sum_{i=1}^d z_i. 
\end{align*}
Fix $\tilde z \in A_n$. Hence from \eqref{ineq:7.4.1}, 
\begin{align*}
\p(\tilde z)& \le \PP(\check Z(\cdot) \ \mbox{first reaches $\tilde C$ before reaching $\tilde A$}| \check Z(0) = \kappa_2 n \tilde z)\\
& \le \PP(\mbox{$L(\check Z(\cdot))$ first reaches $\underline{l}\kappa_2\sum_{i=1}^d z_i$ before reaching $\underline{l}\kappa_2(\epsilon-n\Delta^\eta)$ }| \check Z(0) = \kappa_2 n \tilde z).
\end{align*}
We next consider the following martingale: For $\lambda>0$,
\begin{align*}
U(t) & = \exp\left\{\lambda L( \check Z(t)) - \lambda c(t) - \frac{\lambda^2}{2}\int_0^t|\dot{L}(\check Z(s))|^2 ds \right\},
\end{align*}
where
\begin{align*}
c(t) = \int_0^t \frac{1}{2}tr(\ddot{D}L(\check Z(s))) ds + \int_0^t \sum_{i=1}^d \dot{L}_i(\check Z(s))\theta_i ds + \sum_{i=1}^d \int_0^t \sum_{j=1}^d \dot{L}_j(\check Z(s))R_{ij} dY_i(s).
\end{align*}
Let $\tau(\tilde z) = \inf\{t\ge 0: \kappa_3L(\check Z(t)) \in \tilde A \cup \tilde C\}.$ Now let $m=n=1$ and $\Delta = t$ in \eqref{An} in the proof of Proposition \ref{hitting}, we have for positive constants $\kappa_4$, 
\begin{align*}
c(t) \le - \kappa_4 t, \ 0 \le t \le \tau(\tilde z).
\end{align*}
From Theorem \ref{lyp} (vii), $|\dot{L}(Z(t))| \le M_4$, and we can choose $\lambda$ such that 
\[
\lambda c(t) + \frac{\lambda^2}{2}\int_0^t|\dot{L}(\check Z(s))|^2 ds \le 0, \ 0\le t\le \tau(\tilde z).
\]
Using optional sampling theorem, we have 
\begin{align*}
\exp\{\lambda L(\kappa_2 n \tilde z)\} & = \EE[U(0)] = \EE[U(\tau(\tilde z))] \\
& =  \p(\tilde z) \EE[U(\tau(\tilde z))| L(\check Z(\tau(\tilde z))) = \underline{l}\kappa_2 \sum_{i=1}^d z_i] \\
& \quad + (1-\p(\tilde z)) \EE[U(\tau(\tilde z))|  L(\check Z(\tau(\tilde z))) = \underline{l}\kappa_2 (\epsilon - n \Delta^\eta)]  \\
& \ge  \p(\tilde z) \exp\left\{\lambda \underline{l} \kappa_2 \sum_{i=1}^d z_i \right\} +(1- \p(\tilde z)) \exp\left\{\lambda \underline{l}\kappa_2 (\epsilon - n\Delta^\eta)\right\}.
\end{align*}

Hence
\begin{align}
\p(\tilde z) & \le \frac{\exp\{\lambda L(\kappa_2 n \tilde z)\}- \exp\left\{\lambda  \underline{l} \kappa_2(\epsilon - n\Delta^\eta)\right\} }{\exp\left\{\lambda \underline{l} \kappa_2 \sum_{i=1}^d z_i \right\} -  \exp\left\{\lambda \underline{l} \kappa_2 (\epsilon - n\Delta^\eta)\right\} } \nonumber\\
& \le \frac{\exp\{\lambda \kappa_2 \epsilon L(n\tilde z/\epsilon)\}- \exp\left\{\lambda  L(n\tilde z/\epsilon)\kappa_2(\epsilon - n\Delta^\eta)\right\} }{\exp\left\{\lambda \underline{l} \kappa_2 \sum_{i=1}^d z_i \right\} -  \exp\left\{\lambda \underline{l} \kappa_2 (\epsilon - n\Delta^\eta)\right\} } \nonumber\\
& \le (1- \exp\{-n\Delta^\eta \lambda \underline{l} \kappa_2\}) \frac{\exp\{\lambda \kappa_2 \epsilon L(n\tilde z/\epsilon)\}}{\exp\left\{\lambda\underline{l}  \kappa_2 \sum_{i=1}^d z_i \right\} -  \exp\left\{\lambda\underline{l}  \kappa_2 (\epsilon - n\Delta^\eta)\right\}} \nonumber\\
& = O(n\Delta^\eta). \label{ineq:7.4.2}
\end{align}
Applying \eqref{ineq:7.4.2} and \eqref{ineq:7.4.3} to \eqref{ineq:7.4.4}, we have that
\[
\frac{\max_{\tilde z\in C_n}p_n^-(\tilde z)}{\max_{\tilde z \in C_n} p_n(\tilde z)}  \le 1 + O(n\Delta^\eta).
\]
This completes the proof. \qed

{\bf Proof of Proposition \ref{error}:} We note that 
\begin{align*}
&  \frac{p^\Delta_n(z_n)}{\max_{\tilde z\in C_n}p_n(\tilde z)} - 1 \\
&=\frac{p^\Delta_n(z_n)}{\max_{\tilde z\in C_n}p_n^-(\tilde z)}\frac{\max_{\tilde z\in C_n}p_n^-(\tilde z)}{\max_{\tilde z\in C_n}p_n(\tilde z)}-1
\end{align*}
Recall that $\tau_n^- = \inf\{t\ge 0: \sum_{k=1}^d Z_{n,k}(t) \ge 1-\Delta^\eta, \ \mbox{or} \ \sum_{k=1}^d Z_{n,k}(t) \le \epsilon/n - \Delta^\eta\}$, and $p_n^-(z_n) = \PP(\sum_{k=1}^dZ_{n,k}(\tau_n^-) = 1-\Delta^\eta | Z_n(0) = z_n).$  We note that Lemma \ref{exp_mome2} is still valid when $\tau_n$ is replaced by $\tau_n^-$. Thus, from the proof of Lemma \ref{exp_mome2}, we have 
\[
\PP \left( \sup_{0\le s \le \tau_n-} |Z_n(s)-Z_n^{\Delta}(s)|>{\Delta}^{\eta}\right) = e^{-n \xi^2/(2\Delta) + o(n)},
\]
where $\xi = \frac{{{\Delta}^{\eta}}/C - \Delta|\theta|}{d|\Sigma^{1/2}|}.$
Also Theorem \ref{SPLDP} is valid for $p_n^-(z_n)$, which yields that 
\[
p_n^-(z_n) = e^{-n \inf_{x\in B} I(x) + o(n)}.
\]

We observe that 
\begin{align*}
p_n^{\Delta}(z_n)& = \PP\left(\tau_n^{\Delta}(B_n)<\tau_n^{\Delta}(A_n),\sup_{0\le s \le \tau_n^-}|Z_n(s) - Z_n^\Delta(s)| < {\Delta}^{\eta}\right)\\
& \quad \quad +\PP\left(\tau_n^{\Delta}(B_n)<\tau_n^{\Delta}(A_n),\sup_{0\le s \le \tau_n^-}|Z_n(s) - Z_n^\Delta(s)| > {\Delta}^{\eta}\right)\\
& \le p_n^{-}(z_n)+\PP\left(\sup_{0\le s \le \tau_n^-}|Z_n(s) - Z_n^\Delta(s)| > {\Delta}^{\eta}\right)
\end{align*}
Noting that $\frac{\xi^2}{2\Delta}=\Theta(\Delta^{2\eta-1}),$ for some constant $C>0$, we have
\begin{align*}
\frac{p_n^{\Delta}(z_n)}{p_n^-(z_n)}\le 1+ \frac{\PP\left(\sup_{0\le s \le \tau_n^-}|Z_n(s) - Z_n^\Delta(s)| > {\Delta}^{\eta}\right)}{p_n^-(z_n)} =1 + o\left(e^{-Cn\Delta^{2\eta-1}}\right).
\end{align*}
At last from Lemma \ref{san2}, we obtain 
\[
 \frac{\max_{\tilde z\in C_n}p^-_n(\tilde{z})}{\max_{\tilde z\in C_n}p_n(\tilde z)} \le  1 + O\left(n\Delta^{\eta}\right).
\] 
Combining the above estimates, we establish that 
\[\frac{\max_{\tilde z\in C_n} p^\Delta_n(\tilde z)}{\max_{\tilde z\in C_n}p_n(\tilde z)} \le 1 +  O\left(n\Delta^{\eta}\right).\]
To show the other direction, we define $\tau_n^+ = \inf\{t\ge 0: \sum_{k=1}^d Z_{n,k}(t) \ge 1+\Delta^\eta, \ \mbox{or} \ \sum_{k=1}^d Z_{n,k}(t) \le \epsilon/n + \Delta^\eta\}$, and $p_n^+(z_n) = \PP(\sum_{k=1}^dZ_{n,k}(\tau_n^+) = 1+\Delta^\eta | Z_n(0) = z_n).$ Similar to Lemma \ref{san2} (replace $p^-_n(z_n)$ and $p_n(z_n)$ by $p_n(z_n)$ and $p^+_n(z_n)$, respectively, in the proof), it can be shown that 
\[
 \frac{\max_{\tilde z\in C_n}p^+_n(\tilde{z})}{\max_{\tilde z\in C_n}p_n(\tilde z)} \ge  1 + O\left(n\Delta^{\eta}\right).
\]
Further note that 
\begin{align*}
p_n^{\Delta}(z_n)& = \PP\left(\tau_n^{\Delta}(B_n)<\tau_n^{\Delta}(A_n),\sup_{0\le s \le \tau_n^-}|Z_n(s) - Z_n^\Delta(s)| < {\Delta}^{\eta}\right)\\
& \quad \quad +\PP\left(\tau_n^{\Delta}(B_n)<\tau_n^{\Delta}(A_n),\sup_{0\le s \le \tau_n^-}|Z_n(s) - Z_n^\Delta(s)| > {\Delta}^{\eta}\right)\\
& \ge p_n^{+}(z_n).
\end{align*}
Finally, we have 
\begin{align*}
&  \frac{\max_{\tilde z\in C_n} p^\Delta_n(\tilde z)}{\max_{\tilde z\in C_n}p_n(\tilde z)}  =\frac{\max_{\tilde z\in C_n} p^\Delta_n(\tilde z)}{\max_{\tilde z\in C_n}p_n^+(\tilde z)} \ \frac{\max_{\tilde z\in C_n}p_n^+(\tilde z)}{\max_{\tilde z\in C_n}p_n(\tilde z)} \ge 1 +  O\left(n\Delta^{\eta}\right).
\end{align*}

\qed

\subsubsection{Proofs of Propositions \ref{alg_stable} and \ref{restart_stable}}

{\bf Proof of Theorem \ref{alg_stable}.} The stability follows directly from Proposition 7 in \cite{dean2009}, and \eqref{split_unbiased} follows from Lemma 5 in the same reference.  To prove \eqref{sp_var}, we need to verified that Conditions 1 and 4 in \cite{dean2009} hold for our problem. Condition 1 of  \cite{dean2009} follows from Theorem \ref{SPLDP}. For their Condition 4, due to the positive recurrence of our defined process, $\sigma_n$ is finite almost surely. Then by strong Markov property

\begin{align*}
 1_{\{Z_n(\sigma_n)\in L_x\}}(p_n(Z_n(\sigma_n)))^2\le 1_{\{Z_n(\sigma_n)\in L_x\}}(p_n(Z_n(\sigma_n))) \sup\limits_{y\in \partial L_x}p^n(y).
\end{align*} 
Hence using the strong Markov property of $Z_n$, we have 
\begin{align*}
\EE_{z_n}\left[1_{\{Z_n(\sigma_n)\in L_x\}}(p_n(Z_n(\sigma_n)))^2 \right]&\le \sup\limits_{y\in \partial L_x}p_n(y) \cdot \EE_{z_n}\left[ 1_{\{Z_n(\sigma_n)\in L_x\}}p_n(Z_n(\sigma_n))\right]\\
&=\sup\limits_{y\in \partial L_x}p^n(y) \cdot p_n(z_n).
\end{align*}
By the large deviation upper bound, we can get 
\begin{align*}
& \liminf\limits_{n\rightarrow \infty}-\frac{1}{n}\log \EE_{z_n}\left[ 1_{\{Z_n(\sigma_n)\in L_x\}}(p_n(Z_n(\sigma_n)))^2 \right]\\
&=-\limsup\limits_{n\rightarrow \infty}\frac{1}{n}\log \EE_{z_n}\left[ 1_{\{Z_n(\sigma_n)\in L_x\}}(p_n(Z_n(\sigma_n)))^2 \right]\\
&\ge -\limsup\limits_{n\rightarrow \infty}\frac{1}{n}\log\left[ \sup\limits_{y\in \partial L_x}p_n(y) \cdot p_n(z_n) \right]\\
&=  -\limsup\limits_{n\rightarrow \infty}\frac{1}{n}\left[ \log \sup\limits_{y\in \partial L_x}p_n(y) +\log p_n(z_n) \right]\\
&\ge -\limsup\limits_{n\rightarrow \infty}\frac{1}{n} \log \sup\limits_{y\in \partial L_x}p_n(y) -\limsup\limits_{n\rightarrow \infty}\frac{1}{n}\log p_n(z_n) \\
&\ge \inf\limits_{y\in \partial L_x}W(y) + W(z).
\end{align*}
Now \eqref{sp_var} follows directly from Theorem 8 in \cite{dean2009}. \qed
\\

\noindent {\bf Proof of Proposition \ref{restart_stable}:} 
It suffices to verify the five conditions in \cite {dean2011}, which are Conditions 3.1, 4.3, 5.1, 5.2, 5.3, are satisfied. We note that for our specific problem, those conditions essentially require that:
\begin{itemize}
	\item $A_n$ and $B$ are closed sets.
	\item $Z_n^{\Delta}(t)$ satisfies the large deviation principle.
	\item For any compact $K\subset \RR^d_+$, all zero cost trajectories with initial positions in $K$ will enter $(A_n\cup B)^{\circ}$ by some fixed finite time depended on $K$.
\end{itemize}
Indeed all these conditions are satisfied, and thus the proposition follows.
\qed

\bibliographystyle{plain}
\bibliography{ref}

\vspace{0.5in}

\noindent Kevin Leder and Zicheng Wang \\
Department of industrial and system engineering\\
University of Minnesota\\
Minneapolis, MN 55455, USA \\
kevin.leder@isye.umn.edu, wang2569@umn.edu.

\vspace{0.2in}

\noindent Xin Liu\\
School of Mathematicla and Statistical Sciences\\
Clemson University\\
Clemson, SC 29678, USA\\
xliu9@clemson.edu.

\end{document}